%% file: PSLH-final-v2.tex
\documentclass{article}

\usepackage{amsmath,amssymb,graphicx,longtable}

\usepackage[
      colorlinks=true,
      linkcolor=blue,
      urlcolor=blue,
      filecolor=blue,
      citecolor=red,
      pdfstartview=FitV,
      pdftitle={},
      pdfauthor={},
      pdfsubject={},
      pdfkeywords={},
      pdfpagemode=None,
      bookmarksopen=true
]
{hyperref}

\usepackage{multicol,color}
\definecolor{darkred}{rgb}{0.5,0,0.5}
\definecolor{darkgreen}{rgb}{0,.5,0}
\hypersetup{pdfborder={0 0 0},colorlinks=true,urlcolor=darkred,citecolor=blue,linkcolor=darkgreen}

\newcommand{\be}{\begin{eqnarray}}
\newcommand{\ee}{\end{eqnarray}}
\newcommand{\nn}{\nonumber}

\newcommand{\Oo}{\mathsf{O}}
\newcommand{\Hh}{\mathsf{H}}
\newcommand{\Cc}{\mathsf{C}}

\newcommand{\reals}{\mathbb{R}}
\newcommand{\cx}{\mathbb{C}}
\newcommand{\ints}{\mathbb{Z}}
\newcommand{\oct}{\mathbb{O}}
\newcommand{\quat}{\mathbb{H}}
\newcommand{\Div}{\mathbb{A}}

\newcommand{\Ba}{\bar{a}}
\newcommand{\Bb}{\bar{b}}
\newcommand{\Bc}{\bar{c}}
\newcommand{\Bd}{\bar{d}}

\newcommand{\Bz}{\bar{z}}

\newcommand{\Bu}{\bar{u}}

\newcommand{\ve}{\varepsilon}
\newcommand{\LB}{{\bigtriangleup_{\text{LB}}}}

\newcommand{\mg}{{\mathfrak{g}}}

\newcommand{\la}{\lambda}
\newcommand{\cH}{\mathcal{H}}
\newcommand{\cE}{\mathcal{E}}
\newcommand{\cP}{\mathcal{P}}
\newcommand{\cO}{\mathcal{O}}
\newcommand{\cQ}{\mathcal{Q}}

\newcommand{\cT}{\mathcal{T}}
\newcommand{\capX}{{{X}}}
\newcommand{\capY}{{{Y}}}

\newcommand{\sB}{{\scriptscriptstyle{B}}}

\numberwithin{equation}{section}
\newtheorem{thm}{Theorem}[section]
\newtheorem{lemma}[thm]{Lemma}
\newtheorem{cor}[thm]{Corollary}
\newcommand{\qed}{\hspace{\stretch{1}} $\square$ \\
  \noindent}
\newcommand{\Pf}{\noindent \textbf{Proof. }}

\newcommand{\mf}[1]{\mathfrak{#1}}

\begin{document}

{\flushright 
{ULB-TH/10-26}\\
{AEI-2010-115}\\
}

\vskip .7 cm
\hfill
\vspace{10pt}
\begin{center}
{\large {\bf Modular realizations of hyperbolic Weyl groups}}

\vspace{1.1cm}

{\centering \rule[0.1in]{13cm}{0.5mm} }

\vspace{.5cm}

{\bf Axel Kleinschmidt${}^{1,2}$, Hermann Nicolai\footnotemark[2] and Jakob Palmkvist\footnotemark[1]}

\vspace{10pt}
\footnotemark[1]{\em Physique Th\'eorique et Math\'ematique\\
Universit\'e Libre
de Bruxelles \& International Solvay Institutes\\ 
ULB-Campus Plaine C.P. 231, BE-1050 Bruxelles, Belgium}\\[5mm]

\footnotemark[2]{\em Max-Planck-Insitut f\"ur Gravitationsphysik, Albert-Einstein-Institut\\
Am M\"uhlenberg 1, DE-14476 Potsdam, Germany}\\[10mm]

\end{center}

\vspace{3pt}
\begin{center}
\textbf{Abstract}\\[5mm]
\parbox{13cm}{\footnotesize
We study the recently discovered isomorphisms between hyperbolic Weyl groups 
and modular groups over integer domains in normed division algebras. We show how to realize
the group action via fractional linear transformations on
generalized upper half-planes over the division algebras,
focussing on the cases involving quaternions and octonions. For these we construct automorphic forms, whose explicit expressions
depend crucially on the underlying arithmetic properties 
of the integer domains. Another main new result is the explicit octavian realization of $W^+({\rm E}_{10})$, which contains as a special case a new realization of $W^+({\rm E}_8)$ in terms of unit octavians and their automorphism group.}
\end{center}

\thispagestyle{empty}
\newpage
\setcounter{page}{1}

\tableofcontents


\section{Introduction}

Following the work of Feingold and Frenkel~\cite{FeFr83} it was realized 
in~\cite{Feingold:2008ih} that very generally there are isomorphisms between 
the Weyl groups of hyperbolic Kac-Moody algebras and (finite extensions or 
quotients of) modular groups over integer domains in the division algebras $\Div$ of real, complex, quaternionic or octonionic numbers; $\Div=\reals,\cx,
\quat,\oct$. The simplest case is associated with the rational integers 
$\ints\subset\reals$: the modular group ${\rm PSL}(2,\ints)$ is isomorphic 
to the even part of the Weyl group of the canonical hyperbolic extension 
of the Lie algebra $\mf{sl}(2,\reals)$~\cite{FeFr83}. This is but the first 
example of a rich and interesting family of novel isomorphisms that rely 
on the existence of integers within the division algebras of higher 
dimensions~\cite{Feingold:2008ih}, culminating in the relation 
$W^+({\rm E}_{10})\cong {\rm PSL}(2,\Oo)$ between the even Weyl group of ${\rm E}_{10}$ (the hyperbolic 
extension of the exceptional Lie algebra ${\rm E}_8$) and the integer octonions.

These isomorphisms are interesting for several reasons. The underlying 
arithmetic structure of the various integer domains could help in 
understanding the structure of the still elusive hyperbolic Kac-Moody 
algebras, especially for the `maximally extended' hyperbolic algebra ${\rm E}_{10}$
which is expected to possess very special properties.\footnote{For instance,
  the arithmetic structure of the modular group ${\rm PSL}(2,\Oo)$ 
  may impose more stringent constraints on the characters 
  and root multiplicities of ${\rm E}_{10}$ than the corresponding modular 
  groups do for the lower rank hyperbolic algebras. ${\rm E}_8$ and its hyperbolic
  extension ${\rm E}_{10}$ are also special because they are the only eligible algebras
  with (Euclidean or Lorentzian) even self-dual root lattices.}
This is certainly a tantalizing possibility but we will not pursue it 
further in this work. Instead, we focus on the modular theory associated 
with the novel modular groups that appear in the isomorphisms, and the
associated `generalized upper half planes'. More specifically, we will 
discuss the quaternionic and the octonionic cases, corresponding to 
the hyperbolic extensions of the split real ${\rm D}_4\equiv\mf{so}(4,4)$ 
and ${\rm E}_{8}$ Lie algebras. In the quaternionic case in particular, we 
give a detailed description of the `generalized upper half plane' of 
dimension five on which the modular group acts naturally. 

Our main motivation derives from potential applications of this 
modular theory in quantum gravity and M-theory. As shown recently,
the generalized upper half planes are the configuration spaces for 
certain models of (anisotropic) mini-superspace quantum gravity, and 
a subclass of the automorphic forms associated with the modular 
groups studied here appear as solutions to the Wheeler-DeWitt 
equation in the  cosmological billiards 
limit~\cite{Kleinschmidt:2009cv,Kleinschmidt:2009hv,Forte}, 
a version of mini-superspace quantum gravity which we refer to as 
`arithmetic quantum gravity'. Automorphic functions 
(particularly Eisenstein series) also feature prominently in recent studies 
of non-perturbative effects in string and M-theory, 
see \cite{Green:1997tv,Obers:1999um} for early work and 
\cite{Green,Lambert:2010pj,Pioline} for recent progress concerning the split
exceptional groups $G={\rm E}_{7}$ and ${\rm E}_{8}$, with 
corresponding arithmetic groups ${\rm E}_7(\ints)$ and 
${\rm {\rm E}_8}(\ints)$.\footnote{In some cases, the relevant `automorphic 
forms' are not eigenfunctions of the Laplacian but have 
source terms~\cite{Green:2005ba}.} 
There exist far
reaching conjectures concerning the infinite-dimensional extensions 
${\rm E}_9$, ${\rm E}_{10}$ and ${\rm E}_{11}$, but it is less clear how to make sense of (or even 
{\em define}) their discrete subgroups (see, however,~\cite{Green:2010kv}). Nevertheless, our results 
can be viewed as a first step towards extending these ideas to the 
infinite-dimensional duality group ${\rm E}_{10}$, since ${\rm PSL}(2,\Oo)$, being the even Weyl group, 
is expected to be contained in any hypothetical arithmetic 
subgroup of ${\rm E}_{10}$. A possible link with the work of 
\cite{Kleinschmidt:2009cv,Kleinschmidt:2009hv,Forte} 
is provided by a conjecture of \cite{Ganor} according to which 
the solution of the Wheeler-DeWitt equation in M-theory is a 
vastly generalized automorphic form with respect to ${\rm E}_{10}(\ints)$.
Our results on ${\rm D}_4\equiv\mf{so}(4,4)$
can also be viewed 
from a string theory perspective, since ${\rm SO}(4,4)$ can be realized as the symmetry of type I supergravity in ten dimensions without any vector multiplets on a four-torus. It is also the continuous version of T-duality of type II strings.

Within the general theory of automorphic forms for groups of real rank one, the main new feature that distinguishes our 
construction  is  the link with the division algebras $\Div$  and the
integer domains $\cO\subset\Div$ (which always contain $\ints$ 
as a `real' subset). Moreover, the arithmetic modular groups considered  
here are all identified with the even subgroups of certain hyperbolic 
Weyl groups. That is, each of the modular groups considered here is 
isomorphic to the {\em even} Weyl group of a canonical hyperbolic 
extension $\mg^{++}$ of an associated simple finite-dimensional 
split Lie algebra $\mg$,
\begin{align}\label{Gamma}
\Gamma \equiv W^+_{\text{hyp}}\equiv W^+(\mg^{++}).
\end{align}
More specifically, the groups $\Gamma$ are discrete subgroups 
of the groups ${\rm SO}(1,n\!+\!1)$ acting on the (Lorentzian)
root space of $\mg^{++}$, but with additional restrictions implying 
special  (and as yet mostly unexplored) arithmetic properties. 
The standard Poincar\'e series construction then implies a simple 
general expression for an automorphic function on the generalized 
upper half plane as a sum over images of the group action as
\begin{align}
f(u,v) = \sum_{\gamma\in \Gamma_\infty\backslash \Gamma} (\gamma\cdot v)^s,
\end{align}
where $u\in\Div$ and $v\in \reals_{>0}$ parametrize the upper half plane.
(This expression is convergent for $\text{Re}(s)$ sufficiently large 
and can be analytically continued to most complex values of $s$ by 
means of a suitable functional relation.) The subgroup $\Gamma_\infty$ 
appearing in this sum is conventionally defined as the subgroup of the modular
group $\Gamma$ leaving invariant the cusp at infinity; here this group
becomes identified with the even Weyl group of the affine subalgebra
$\mg^+\subset \mg^{++}$,
\begin{align}\label{Gamma1}
\Gamma_\infty \equiv W_{\text{aff}}^+\equiv W^+(\mg^+) \subset \Gamma.
\end{align}
As we will explain, the cusp at infinity in the generalized upper half 
plane here corresponds to the affine null root of the affine algebra
$\mg^+$, in accord with the fact that $W_{\text{aff}}$ can be defined 
as the subgroup of $W_{\text{hyp}}$ stabilizing the affine null root~\cite{GN}.
At the same time $W_{\text{aff}}^+$ is the {\em maximal parabolic 
subgroup} of $\Gamma= W_{\text{hyp}}^+$.  Because $W_{\text{aff}} = 
W_{\text{fin}}\ltimes \cT$, where $W_{\text{fin}}\equiv W(\mg)$ and $\cT$ is the
abelian group of affine translations, the minimal parabolic 
subgroup is $\cT\subset\Gamma_\infty$.

We briefly position our work relative to the existing mathematical literature 
that we are aware of. Generalized upper half planes have appeared 
for example in \cite{Goldfeld} where they are defined as cosets 
${\rm GL}(n,\reals)/{\rm O}(n)\times {\rm GL}(1,\reals)$, on which 
the discrete groups ${\rm GL}(n,\ints)$ act as generalized modular groups. 
This definition fits with our definition only for $n=2$ (corresponding 
to $\Div=\reals$). More generally, one can define generalized upper half
planes as quotient spaces $G/K(G)$ where $G$ is a non-compact group
and $K(G)$ its compact subgroup, and then consider the action of some
arithmetic subgroup $G_\ints\subset G$ on this space.  Our discrete 
groups $\Gamma$ are subgroups of ${\rm SO}(1,n+1;\reals)$ and live on the 
well-known symmetric space ${\rm SO}(1,n+1)/{\rm SO}(n+1)$~\cite{Helgason}. 
Therefore our analysis is in principle part of the general theory of reductive groups 
of real rank one (see for example~\cite{Boulder,HC,Langlands,Borel,Iwaniec,ParkCity}). 
What makes it stand out and interesting for us -- apart from the fascinating
potential applications in fundamental physics -- is the precise nature of the discrete 
subgroups $\Gamma$, which are different from the arithmetic groups $G_\ints$
usually considered in this context,\footnote{These are obtained by choosing a metric to define the matrix group ${\rm SO}(1,n+1)$ and then restricting the matrix entries to (rational) integers.} and the  link to integers in normed division 
algebras and to hyperbolic Weyl groups, which has not been exposed in the 
existing literature to the best of our knowledge. For example, it would obscure 
the underlying arithmetic  structure, if one tried to describe our final expression 
(\ref{eisen}) in terms of a complicated lattice sum in $\reals^{2n}$. It is not even 
clear to us how to derive such a sum without the knowledge of the integral structure.

We emphasize that it is not our intention here to give a complete
account of the theory of automorphic forms for the modular groups we study
(after all, even for $\Gamma= {\rm SL}(2,\ints)$ this is a subject which 
easily fills a whole book \cite{Borel,Iwaniec}). Rather, we regard the
present work merely as a first step towards a theory that remains to 
be fully developed in future work. In particular, since we are using a direct 
and specific  construction in terms of Poincar\'e series, we have not explored 
the adelic approach to modular forms (see for example~\cite{Bump}) although 
we anticipate that this could yield interesting results in this context. 

This article is structured as follows. In section~\ref{sec:uhp}, we review 
the construction of normed division algebras and introduce the generalized 
upper half planes that we will utilize in the remainder of the text. 
Section~\ref{sec:modgroups} is devoted to the isomorphism between modular 
groups and hyperbolic Weyl groups in general; section~\ref{sec:d4sec} then 
deals in detail with the case of quaternions and ${\rm D}_4\equiv\mf{so}(4,4)$, whereas
in section~\ref{octsec}, we present the modular group ${\rm PSL}(2,\Oo)$ over 
the non-associative integer octonions $\Oo$ that appears for ${\rm E}_{10}$, 
the hyperbolic extension of ${\rm E}_8$.  In particular, it gives
a novel realization of the even (and finite) Weyl  group $W^+({\rm E}_8)$ 
in terms of the 240 octavian units and the finite automorphism group 
${\rm G}_2(2)$ of the octavians. Automorphic forms for all these groups 
are defined and analysed in detail in section~\ref{sec:quatauto}, with 
special emphasis on the most interesting case $\Div = \oct$. These sections 
form the heart of this paper and contain the bulk of our new results. 
Several appendices contain some additional results and technical details that we use in the 
main text.

\section{Division algebras and upper half planes}
\label{sec:uhp}

\subsection{Cayley-Dickson doubling}
\label{CDsub}

The division algebras
of complex numbers, quaternions and octonions can be defined 
recursively from the real numbers by a procedure called 
{\it Cayley-Dickson doubling} \cite{Baez,Conway}.
Let $\Div$ be a real finite-dimensional algebra with a conjugation, 
i.e.~a vector space involution $x \mapsto \bar x$ such that 
$\overline{ab}=\bar b \, \bar a$ for any $a,\,b \in \Div$.
In the Cayley-Dickson doubling
one introduces a vector space $i\Div$, isomorphic to $\Div$, 
where $i$ is a new imaginary unit,  and considers the direct sum 
$\Div \oplus i\Div$ of these two vector spaces. The product of 
two elements in $\Div \oplus i\Div$ is defined as \cite{Conway}
 \begin{align} \label{doublingformula-mult}
 (a+ib) (c+i d) &:= (ac -d\bar{b}) +i (cb+\bar{a}d),
 \end{align}
and conjugation as
\begin{align} \label{doublingformula-conj}
\overline{a+i b}:=\bar{a}-i b
\quad\qquad \mbox{(with $ai=i\Ba$).}
\end{align}
Starting from 
$\Div=\reals$ (with the identity map as conjugation) and successively applying the Cayley-Dickson doubling one thus obtains an infinite sequence of so-called
{\it Cayley-Dickson algebras}. 
The first 
doubled algebra is of course the familiar algebra $\cx$ of complex numbers.
The following two algebras in the sequence are $\quat$ and $\oct$, consisting of quaternions and octonions, respectively.
Studying these algebras one finds that $\quat$ is associative but not commutative, whereas $\oct$ is neither associative nor commutative. Nevertheless, 
$\oct$ is an {\it alternative} algebra, which means that any subalgebra generated by two elements is associative.

Any Cayley-Dickson algebra admits a positive-definite inner product
\begin{align}\label{ab}
(a,b):=\frac12 (a\bar b +b \bar a)
\end{align}
with the associated norm $|a|^2\equiv (a,a) = a\bar a$. As a consequence of 
(\ref{doublingformula-mult}) and (\ref{doublingformula-conj}) the norm of the `doubled' expression 
$a+ib\in \Div \oplus i\Div$ is 
\begin{align}\label{CD1}
|a +ib|^2 = |a|^2 + |b|^2.
\end{align}
From the doubling formulas (\ref{doublingformula-mult}) and (\ref{doublingformula-conj}) it follows
\begin{align}
|(a+ib)(c+id)|^2 &= |ac-d\Bb|^2 + |cb + \Ba d|^2 \nn\\
  &\!\!\!\!\!\!\!\!\!\!\!\!\!\!\!\!\!\!\!\!\!\!\!\!\!\!\!\!
  = (|a|^2 + |b|^2)(|c|^2 +|d|^2)    
  +   \,   2\, {\rm Re}  \Big( (\Bd a)(cb) - (b\Bd)(ac)\Big),
\end{align}                 
where
${\rm Re}\, a := \frac12 (a + \Ba)$.
Because for $a,b,c\in \Div$ (for all $\Div$)
\begin{align}
{\rm Re} \,a(bc) &= {\rm Re} \, (ab)c, &
{\rm Re} \, ab &= {\rm Re} \, ba,
\end{align}
we get
\begin{align}\label{CD2}
|(a+ib) (c+id)|^2 &= |a+ib|^2 |c +id|^2 
+ \, 2 \, {\rm Re} \Big( \Bd\, \{a,c,b\}\Big),
\end{align}
where the associator $\{a,b,c\}:=a(bc) - (ab)c$ vanishes for 
associative algebras. The fact that $\reals,\,\mathbb{C}$ and 
$\quat$  are associative thus implies that the composition property
\begin{align}
|ab|=|a||b| \label{compprop}
\end{align}
holds for their doubled algebras $\cx,\,\quat$ and $\oct$. 
Algebras with this property are called {\it normed division algebras}.
It follows from (\ref{compprop}) that 
if 
$ab=0$ then
$a$ or $b$ must be zero,
so these algebras have no zero divisors.
Hurwitz' theorem (see \cite{Baez,Conway}) states that $\reals$, $\cx$, $\quat$ and $\oct$ 
are the only (real and finite-dimensional) normed division algebras 
with unity.
For instance, since the octonions are non-associative, their
Cayley-Dickson double does not
satisfy the norm composition property (\ref{compprop}).
This is a real 16-dimensional algebra, often
called the {\it sedenions}. 
It is of course non-associative (since it contains the octonions) but also non-alternative, and does possess zero divisors.

\subsection{Generalized upper half planes}

For every normed division algebra $\Div$ we define the {\em generalized 
upper half plane} 
\begin{align}\label{z}
\cH\equiv\cH(\Div) :=\left\{ z= u + iv  \quad\mbox{with $u\in\Div$ and {\em real} $v>0$}\right\}.
\end{align}
By $i$ we will always denote the new imaginary unit not contained in $\Div$,
while $u\equiv u^0 +\sum u^i e_i$ with the imaginary units $e_i$ in $\Div$.
Hence $\cH(\Div)$ is contained in a hyperplane in the Cayley-Dickson
double $\Div\oplus i\Div$, and of (real) dimension
$\dim_\reals \cH(\Div) = (\dim_\reals\,\Div) +1$. By complex 
conjugation (\ref{doublingformula-conj}) we get $\Bz = \Bu - iv$, and
thus $\Bz$ parametrizes  the corresponding `lower half plane' 
$\overline{\cH}$.  From (\ref{CD2}) we see that the composition 
property
\begin{align}
|z
z'|=|z|
|z'| \qquad z,z'\in\cH(\Div)
\end{align}
continues to hold for {\em all} normed division algebras $\Div=\reals,\cx,\quat,\oct$. 
In particular, it still holds for $\Div=\oct$, even though in this case
$z$ is a sedenion, because the associator in (\ref{CD2}) still vanishes 
as long as $v$ remains real. Likewise, the alternative laws
\begin{align} \label{altlaws}
(aa)z&= a(az), &
(az)a &= a(za), & (za)a&= z(aa)
\end{align}
remain valid for $a \in \oct$ and $z \in \cH(\oct)$
by Cayley-Dickson doubling 
(\ref{doublingformula-mult}). This will be important below when 
we define modular transformations.

The line element in $\cH(\Div)$ is
\begin{align}\label{ds}
ds^2 = \frac{|du|^2 + dv^2}{v^2} 
\end{align}
for $u\in\Div \, (=\reals,\mathbb{C},\quat$ or $\oct)$
and $v>0$, with $|du|^2 \equiv d\Bu \,du$. Using $v$ 
to parametrize the geodesic, the geodesic equation is
(with $u'\equiv du/dv$)
\begin{align}
\frac{d}{dv} \left( \frac1{v \sqrt{1 +|u'|^2}} \frac{du}{dv}\right) = 0.
\end{align}
It is then straightforward to see that geodesics are given by straight 
lines parallel to the `imaginary' ($=v$) axis, or by half circles 
starting at $u_1\in\Div$ and ending at $u_2\in\Div$ on the boundary 
$v=0$ of $\cH(\Div)$. The equation of this geodesic half circle for
$z=u +iv \in \cH(\Div)$ reads
\begin{align}
(u-u_1)(\Bu - \Bu_2) + v^2 = 0.
\end{align}
The length of any geodesic segment connecting two
points $z_1, z_2\in \cH(\Div)$ is  
\begin{align}
d(z_1,z_2) = \log \frac{|z_1 - z_2^*| + |z_1 - z_2|}
{|z_1 - z_2^*| - |z_1 - z_2|}
\end{align}
or, equivalently,
\begin{align}\label{Length}
d(z_1,z_2) = 2 \, {\rm artanh}\,
\frac{|z_1 - z_2|}{|z_1- z_2^*|} = \text{arcosh} \left(1+\frac{|z_1-z_2|^2}{v_1v_2}\right),
\end{align}
where we defined
\begin{align}
z=u+iv \quad \Rightarrow\quad z^*:= u-iv \in \overline{\cH(\Div)}\,.
\end{align}
The formula (\ref{Length}) generalizes the familiar formula from complex 
analysis to all division algebras. Likewise, the volume element on 
$\cH(\Div)$ is given by
\begin{align}\label{Volume}
d{\rm vol}(z) := \frac{d^nu\, dv}{v^{n+1}} 
\end{align}
and the Laplace-Beltrami operator reads\footnote{Note that $\frac{\partial}{\partial u}\frac{\partial}{\partial \bar{u}} |u|^2 \equiv \sum_{i=0}^{n-1}\frac{\partial^2}{\partial u_i^2} |u|^2
= 2n$ in our conventions.}
\begin{align}\label{LB}
\bigtriangleup_{\text{LB}} = v^{n+1} \frac{\partial}{\partial v}\left(
     v^{1-n} \frac{\partial}{\partial v}\right) + 
     v^2\frac{\partial}{\partial u}\frac{\partial}{\partial \Bu},
\end{align}
with $n=1,\,2,\,4$ and $8$ for $\Div=\reals,\,\cx,\,\quat$ and $\oct$, 
respectively. 

As far as the geometry is concerned, the `generalized upper half planes' 
$\cH(\Div)$ are special examples of the general coset spaces
\be
\cH_n = \frac{{\rm SO}(1,n+1)}{{\rm SO}(n+1)}
\ee
for $n=1,2,4,8$. These are all hyperbolic spaces of constant negative
curvature, which can be embedded as unit hyperboloids
\be\label{cHn}
\cH_n = \big\{ x^\mu\in \reals^{1,n+1} \,\big|\,  x^\mu x_\mu = 
              - x^+ x^- + x\!\cdot\! x = -1 \,;\, x^\pm >0 \big\}
\ee
into the forward lightcone of ($n$+2)-dimensional Minkowski space 
$\reals^{1,n+1}$, with $x \in\reals^n$. These hyperboloids are isometric 
to the upper half planes introduced above by means of the mapping
\begin{align}\label{uhpco}
x^- &= \frac1{v}\;,&
x^+ &= v + \frac{|u|^2}{v} \;,&
x&= \frac{u}{v}.
\end{align}
where the last equation identifies the components $x^j = u^j/v$ in
the Euclidean subspace $\reals^n$.
The line element (\ref{ds}) is the pull-back of
the Minkowskian line element $ds^2 = -dx^+ dx^- + dx \!\cdot\! dx$.
Similarly, it is straightforward to check that
\be
\int d{\rm vol} (z) \big( \cdots \big) = \int  d^nx \, dx^+\,   dx^- \,
  \delta\big(x^+ x^- - \bar{x} x -1 \big)\big(\cdots \big)
\ee
and to derive the Laplace-Beltrami operator (\ref{LB}) from the 
Klein-Gordon operator in $\reals^{1,n+1}$. Consequently, the geodesic 
length (\ref{Length}), the volume element (\ref{Volume}) and 
the Laplace-Beltrami 
operator (\ref{LB}) are left invariant under the isometry group ${\rm SO}(1,n+1)$ 
of the embedding space (and the unit hyperboloid). Because the even
Weyl groups (or modular groups) to be considered below are all discrete 
subgroups of ${\rm SO}(1,n+1)$, these geometric objects are {\em a fortiori} 
invariant under these discrete groups as well. The main new feature 
here distinguishing the cases $n=1,2,4,8$ from the general case (\ref{cHn})
is the link with the division algebras and their algebraic structure,
which is evident in particular in the form of the modular transformations
(\ref{sz}) below.

For completeness, we discuss the Green function on $\cH(\Div)$ in appendix~\ref{app:Green} and (periodic) geodesics in appendix~\ref{app:geo} for the quaternionic case.

\section{Weyl groups as modular groups}
\label{sec:modgroups}

As vector spaces with a positive-definite inner product, the normed 
division algebras $\reals, \mathbb{C},\,\quat,\,\oct$ can be 
identified with the root spaces of finite-dimensional Kac-Moody algebras 
of rank $1$, $2$, $4$, $8$, respectively. The root lattices are 
then identified with lattices in these algebras. In several cases 
of interest, the lattices close under multiplication and thus endow
the root lattice with the structure of a (possibly non-associative) 
ring~\cite{Feingold:2008ih}.

The fundamental Weyl reflection with respect to a simple root $a$
is defined as a reflection in the hyperplane orthogonal to $a$,
\begin{align}
x \mapsto x-2\frac{(a,\,x)}{(a,\,a)}a\,.
\end{align}
When we identify the roots with elements in the division algebra 
$\Div$ this can be written
\begin{align} \label{fundweylrefl}
x &\mapsto x - |a|^{-2}({a \bar x+x \bar a})a
= x - |a|^{-2}(a \bar x a + x |a|^2)
= -  
\frac{a \bar x a }{|a|^2}\,,
\end{align}
where we have used the definition of the inner product, and the
alternativity of  $\Div$. For simply-laced algebras (the only case
we will consider in this paper), we normalize the simple 
roots to unit norm. In particular, if $\ve_i$ denote the so
normalized simple roots of the finite algebra, then the 
associated finite Weyl group $W_{\text{fin}}$ is generated by the fundamental 
Weyl reflections
\begin{align}\label{wx}
x \mapsto w_i(x)=-\ve_i\bar{x}\ve_i\,,
\end{align}
for $x\in\Div$
with $i$ ranging over the rank of the algebra. We will always choose 
the orientation of the simple roots such that the highest root 
$\theta$ of the algebra is equal to the real unit, i.e., $\theta=1$.

\subsection{Hyperbolic overextensions}

To any simple finite-dimensional Lie algebra (finite Kac-Moody algebra) of rank $r$ one can
associate an infinite-dimensional (indefinite) Kac-Moody algebra of rank $r+2$ 
as follows. 
One first constructs the non-twisted 
affine extension, thereby increasing the rank by one. Then one adds an 
additional node with a single line to the affine node in the Dynkin diagram. 
The resulting algebras 
are often called `over-extended'~\cite{Gaberdiel:2002db,Damour:2002fz} 
and in many cases turn out to be hyperbolic Kac-Moody 
algebras.\footnote{But not always: For instance, the finite algebra 
 ${\rm A}_8$ extends to $A_8^{++}\equiv {\rm AE}_{10}$ 
  which is indefinite, but not hyperbolic.} The Cartan-Killing metric 
on the Cartan subalgebra is always Lorentzian. We will number the 
affine and over-extended nodes by $0$ and $-1$, respectively. Denoting the finite-dimensional split real algebra by $\mf{g}$, the associated over-extension will be denoted by $\mf{g}^{++}$ and we use the index 
$I=-1,\,0,\,1,\,\ldots,\,r$ to denote the simple roots of $\mf{g}^{++}$. The affine extension will be denoted by $\mf{g}^+$.

The infinite Weyl groups associated with these over-extended 
algebras were studied in~\cite{Feingold:2008ih}. The Weyl group 
acts on the Lorentzian vector space $\Div\oplus \reals^{1,1}$, 
{\em alias} the Minkowski space $\reals^{1,n+1}$, consisting of real 
linear combinations of the simple roots. We can identify this vector 
space with the Jordan algebra\footnote{A Jordan algebra 
   is a commutative (but possibly non-associative) algebra 
   where any two elements $X$ and $Y$ satisfy
   $X^2 \circ (Y \circ X)= (X^2 \circ Y) \circ X$.
   This identity holds for $H_2(\Div)$ with $\circ$ being the symmetrized 
   matrix product. However, we will here only consider the Jordan 
   algebras $H_2(\Div)$ as vector spaces, without using the Jordan algebra 
   property.} $H_2(\Div)$ of Hermitean matrices
\begin{align}\label{capXform}
\capX :=  \left( \begin{array}{cc} x^+ & x \\
                                   \bar x & x^- \end{array}\right)
     = \capX^\dagger
\end{align}
where $x^\pm \in \reals$ and $x\in\Div$.
This vector space is Lorentzian with respect to the norm
\begin{align}
\label{nrmq}
||\capX ||^2 := - \det\capX = - x^+ x^- + \bar x x
\end{align}
and the associated bilinear form
\begin{align}
(\capX,\,\capY)=\tfrac12\left(||\capX+\capY||^2-||\capX||^2-||\capY||^2\right).
\end{align}
This definition of the norm differs from the one in \cite{Feingold:2008ih}
by a factor of $1/2$.
We also define the simple roots by
\begin{align}\label{simpleroots}
\alpha_{-1}&=
\begin{pmatrix}
1 & 0\\
0 & -1
\end{pmatrix},&
\alpha_{0}&=
\begin{pmatrix}
-1 & -1\\
-1 & 0
\end{pmatrix},&
\alpha_{i}&=
\begin{pmatrix}
0 & \ve_i\\
\bar \ve_i & 0
\end{pmatrix},
\end{align}
where $|\ve_i|=1$, 
so that they have unit length, $||\alpha_I||^2=1$,
instead of length 2, which is the standard normalization for a simply-laced algebra.
The reason for this is that we want the norm to coincide with the standard norm in the division algebra $\Div$, as soon as we restrict the root space to the finite subalgebra. 

The fundamental Weyl reflections with respect to the simple 
roots (\ref{simpleroots}) are given by~\cite{Feingold:2008ih}
\begin{align}\label{Weyl1}
w_I : \capX \mapsto M_I \bar{\capX} M_I^\dagger
\end{align}
where
\begin{align}
M_{-1}&=
\begin{pmatrix}
0 & 1\\
1 & 0
\end{pmatrix},&
M_{0}&=
\begin{pmatrix}
-1 & 1\\
0 & 1
\end{pmatrix},&
M_{i}&=
\begin{pmatrix}
\varepsilon_i & 0\\
0 & -\bar \varepsilon_i
\end{pmatrix},
\end{align}
and we denote the thus generated hyperbolic Weyl group by $W\equiv W_{\text {hyp}}$.

Since Weyl transformations preserve the norm, we can consider 
their action on elements of fixed norm, and in particular on
those elements that lie on the unit hyperboloid $||\capX||^2 = -1$
inside the forward light-cone of the Lorentzian space. As we showed
in the previous section, this unit hyperboloid is isometric to the 
generalized upper half plane $\cH(\Div)$ via the isometric embedding
(\ref{uhpco}). The lightcone $||\capX||^2 =0$ in $H_2(\Div)$ then
corresponds to (a double cover of) the boundary $\partial\cH(\Div)$, 
consisting of the subspace $\Div$ (that is, $v=0$) and the single 
point $z=i\,\infty$. The latter point is called `cusp at infinity', 
and plays a special role because it is associated 
with the affine null root
\begin{align} \label{nullroot}
\delta=
\begin{pmatrix}
 -1 & 0\\
 0 & 0
\end{pmatrix}\,.
\end{align}
To see this more explicitly we note that the lines with constant $u$
parallel 
to the imaginary axis in $\cH(\Div)$ are obtained via (\ref{uhpco}) 
by projection of the null line 
\begin{align}\label{uhpco1}
x^- &=  1 \;,&
x^+ &= t + |u|^2 \;,&
x&= u
\end{align}
onto the unit hyperboloid inside the forward lightcone of 
$\reals^{1,n+1}$. The tangent vector 
of this null line is then
identified with (\ref{nullroot}) again via (\ref{uhpco}), and
the boundary point on the unit hyperboloid corresponding 
to the cusp $z=i\infty$ is reached for $t\rightarrow\infty$. 
Choosing any other null direction inside the forward 
lightcone one can reach all other points in the boundary,
that is $\Div$ with $v=0$.

Because the Weyl transformations generated by the repeated
action (\ref{Weyl1}) preserve the norm $||\capX||^2$, they 
leave invariant the unit hyperboloid.
As shown in~\cite{Kleinschmidt:2009cv,Kleinschmidt:2009hv} 
the projection of the fundamental reflections of 
the infinite Weyl group $W$ onto the unit hyperboloid,
and hence the generalized upper half plane, induce
the following modular action on $z\in\cH(\Div)$
\begin{align}\label{modularPGL}
w_{-1}(z) &= \frac{1}{\bar{z}}\,,& w_0(z)&=-\bar{z}+1\,,& 
w_i(z)&= -\ve_i\bar{z}\ve_i\,,
\end{align}
where $\varepsilon_i$ are the (unit) simple roots of the 
underlying finite algebra as above. Furthermore, $1/\Bz = z/|z|^2$.

\subsection{Even Weyl group}

We are here mainly interested in the {\em even} part of the Weyl 
group, which is realized by `holomorphic' transformations (see 
(\ref{sz}) below) and which we wish to interpret as a generalized 
modular group. The even Weyl group $W^+\subset W$ consists 
of the words of even length in $W$.  As the generators of the even Weyl 
group  we take the transformations
\begin{align} \label{newevengenerators}
s_I := w_I w_\theta
\end{align}
with $I=-1,0$ or $I=i$, where $w_\theta$ is the reflection on 
the highest root $\theta=1$. 
These generators
act as 
\begin{align}
s_I : \capX \, \mapsto \, S_I {\capX} S_I^\dagger,
\end{align}
where
\begin{align}
\label{evenweylgen}
S_{-1}&=
\begin{pmatrix}
0 & -1\\
1 & 0
\end{pmatrix},&
S_{0}&=
\begin{pmatrix}
1 & 1\\
0 & 1
\end{pmatrix},&
S_{i}&=
\begin{pmatrix}
\varepsilon_i & 0\\
0 & \bar \varepsilon_i
\end{pmatrix}.
\end{align}
On the upper half plane it follows from (\ref{wx}) (with $\ve = \theta =1$)  
that
\begin{align}\label {wtheta}
w_\theta(z) = - \Bz \;,
\end{align}
and the generators (\ref{newevengenerators}) act as
\begin{align}\label{sz}\boxed{
s_{-1}(z) = -\frac{1}{z}, \qquad s_0(z) = z+1, \qquad s_i(z)= \ve_i z \ve_i.}
\end{align}
These transformations are `holomorphic' on the upper half plane 
and generalize the well-known expression for modular transformations 
of the two-dimensional upper half plane $\cH(\reals)$ under 
${\rm PSL}(2,\ints)$.

In \cite{Feingold:2008ih} $\tilde s_0 = w_{-1}w_0$ and $\tilde s_i =w_{-1}w_i$ were used as generating elements. 
Our choice here is more convenient since it refers to the universal element $\theta=1$, but 
unlike the one in \cite{Feingold:2008ih}, the generating set (\ref{newevengenerators}) is not minimal.
For $\Div=\reals, \cx, \quat$ this redundancy is expressed 
in the extra relations 
\be\label{srelations}
s_1 =1 \; \mbox{for $\Div=\reals$}\, , \;\;\;
s_1 s_2 = 1 \; \mbox{for $ \Div = \cx$}\, ,\;\;\;
s_1 s_3 s_4 =1 \;\; \mbox{for $\Div = \quat$}
\ee
with appropriate numbering of the roots, and always in the basis
where $\theta~=~1$. For the octonions the relevant relation is more
tricky because of non-associativity.

Quite generally, the even Weyl groups of the over-extended algebras 
considered in \cite{Feingold:2008ih} are of the form 
\be
W^+\equiv W^+_{\text{hyp}} \cong {\rm PSL}(2,\cO)\,,
\ee 
where $\cO\subset\Div$ is a suitable set of algebraic integers 
in $\Div$. The roots of the associated finite Kac-Moody algebra then 
correspond to the units of $\cO$ (a unit being an invertible element 
of the ring). In some cases one has to consider 
finite extensions of this group or finite quotients (depending on 
whether the algebra is not simply-laced or has diagram automorphisms); 
a detailed analysis of all cases can be found in~\cite{Feingold:2008ih}. 
The even Weyl group $W^+$ is a normal subgroup of index 
two in the full Weyl group and for many purposes it is sufficient to study it.
In the following sections, we will focus on two distinguished cases, one 
related to ${\rm D}_4$ and the quaternions (where $\cO = \Hh$,
the ring of `Hurwitz numbers') and one related to ${\rm E}_8$ and 
the octonions (where $\cO=\Oo$, the `octavians').

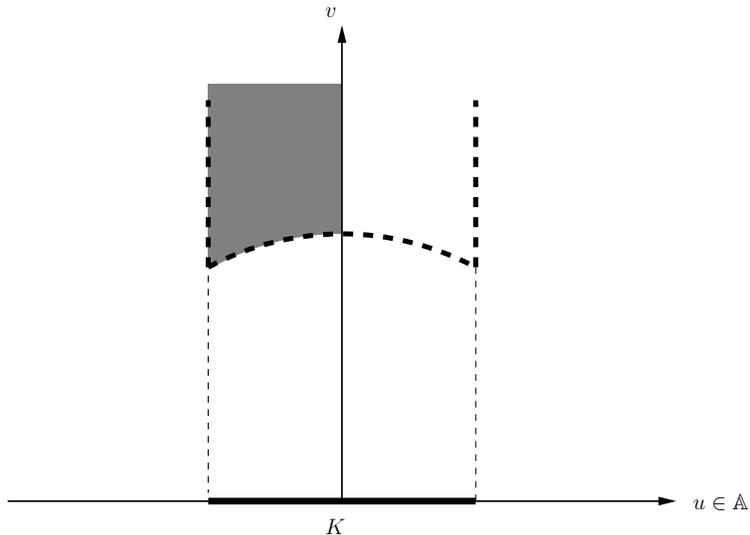
\begin{figure}[h]
\centering
\scalebox{.7}{\input{uhp5.pstex_t}}
\caption{\label{fundfig}\it Sketch of the fundamental domains ${\cal{F}}_0$ 
 (in grey) and $\cal{F}$ (within heavy dashed lines). The shape of $\cal{F}$ 
 is always that of a `skyscraper' over a compact domain $K$ 
 contained in the unit ball of $\Div$. This `skyscraper',  whose bottom 
 has been cut off, extends to infinite height in $v$. }
\end{figure}

The Weyl transformations (\ref{modularPGL}) define a fundamental domain 
${\cal{F}}_0\subset\cH(\Div)$, which is the image of the fundamental Weyl 
chamber in $H_2(\Div)$ under the projection (\ref{uhpco}). Likewise,
the {\em even} Weyl transformations (\ref{sz}) define a fundamental 
domain $\cal{F}$, which contains two copies of ${\cal{F}}_0$, and a 
`projection' $K$ of $\cal{F}$ onto $\Div$; 
see Figure~\ref{fundfig} below for an illustration.  As we 
already explained, the cusp $z=i\infty$ is the image of the 
boundary point of the unit hyperboloid which is reached by 
following and projecting any null ray inside the forward lightcone
along the affine null root $\delta$.

From (\ref{Volume}) it follows that 
\begin{align}
\int_{\cal{F}} d{\rm vol} (z) < \infty
\end{align} 
if $K$ does not extend beyond the unit sphere in the $\Div$-plane. It may touch the unit sphere at isolated points since the improper integral will still be convergent and finite. This happens 
among the over-extended algebras in particular for the algebras
$A_7^{++}, B_8^{++}$ and $D_8^{++}$.\footnote{See~\cite{Fleig:2011mu} for further discussions of the fundamental domain of hyperbolic Kac-Moody algebras and explicit volume computations.}
The finite volume of the fundamental domain is in one-to-one correspondence 
with the hyperbolicity of the underlying Kac-Moody algebra~\cite{Damour:2002fz}.

For all $\cH(\Div)$, the geodesic length (\ref{Length}) is invariant 
under the continuous isometry group ${\rm SO}(1,n+1)$ and therefore also
under the generating elements of the Weyl group (\ref{modularPGL}), 
that is under the full Weyl group, as can also be verified explicitly
for the generating Weyl reflections (\ref{modularPGL}).
For instance, for the inversion we get
\begin{align}
w_{-1}  \; : \quad |z_1 - z_2| \rightarrow |\Bz_2^{-1}-\Bz_1^{-1}| =   
       (|z_1||z_2|)^{-1} |z_1-z_2|  
\end{align}
using alternativity. By iteration it then follows that 
\begin{align}
d(z_1,z_2) = d(w(z_1), w(z_2))
\end{align}
for all elements $w$ of the Weyl group $w\in W_{\text{hyp}}$. 
Similarly, for all $w\in W_{\text{hyp}}$ the hyperbolic volume 
element is left invariant:
\be
d{\rm vol} (z) = d{\rm vol} (w(z)) \,.
\ee
Finally, the Laplace-Beltrami operator is also invariant in the
following sense. Setting $z':=w(z)$ and denoting by 
$\bigtriangleup'_{\text{LB}}$ the Laplace-Beltrami operator with respect to 
the primed coordinates $z'$ we have {\em for any function} $f(z)$ 
and for all $w\in W_{\text{hyp}}$
\be\label{LB1}
\bigtriangleup'_{\text{LB}} f(z') =
\bigtriangleup_{\text{LB}} f(w(z)) \; .
\ee

\subsection{Affine Weyl group}
\label{sec:affweyl}

The hyperbolic Weyl group has a natural subgroup associated with the 
affine subalgebra, corresponding to the embedding $\mf{g}^{+}
\subset\mf{g}^{++}$. We denote this subgroup by $W_{\text{aff}}\equiv
W(\mf{g}^+)$. It is known~\cite{Ka90} that it has the structure
\be\label{affgp}
W_{\text{aff}} = W_{\text{fin}} \ltimes \cO,
\ee
where $\cO$ is the integer domain corresponding to the root lattice of the finite-dimensional algebra $\mf{g}$. In (\ref{affgp}), $\cO$ is to be thought of as a free abelian group $\mathcal{T}$ of translations on the finite root lattice on which the finite Weyl group $W_{\text{fin}}$ acts.
Written in matrix form the translations take the form
\be
T_y = \begin{pmatrix}
1 & y\\
0 & 1
\end{pmatrix},
\ee
and obey
$T_x T_y = T_{x+y}$ for all $x,y\in\cO$. In the associative cases one can write the finite Weyl transformations as matrices as
\be
M_{a,b}=\begin{pmatrix}
a & 0\\
0 & b
\end{pmatrix}
\ee
for units $a,b\in\cO$, possibly with restrictions depending on $\mf{g}$.\footnote{We will encounter an explicit example of such restrictions in the case $\mf{g}={\rm D}_4$ below.} It is easy to see that the transformations leaving invariant the null root (\ref{nullroot}) are exactly the $T_x$ and the $M_{a,b}$. This is to be expected as the affine Weyl group can be defined as the stabilizer of the affine null root~\cite{GN}.

\section{${\rm D}_4$ and its hyperbolic extension} \label{quat-section}
\label{sec:d4sec}

Now we specialize to $\Div=\quat$ in which case the
relevant ring of algebraic integers is $\cO = \Hh$,
the `Hurwitz numbers', to be defined below. The results
of this section generalize many well known results
for the usual upper half plane ($\Div=\reals$). Since they
contain similar results for $\Div=\cx$ we will not specially consider
this case (which has been dealt with in the literature \cite{Then}).
Our main task will be to find closed formulas for arbitrary elements
of the group ${\rm PSL}(2,\Hh)$ and its modular realization
on $\cH(\quat)$. In comparison to the commutative cases 
$\Div=\reals,\cx$ this case presents new features
because of the non-commutativity of the matrix entries.
However, a direct generalization to $\Div=\oct$ is not possible 
because of the non-associativity of the matrix entries; this case 
will be treated separately in section~\ref{octsec}.

\subsection{${\rm D}_4$ and the integer quaternions}

The quaternions are obtained by applying the Cayley-Dickson doubling 
procedure described in section~\ref{CDsub} to $\Div=\cx$. An arbitrary 
quaternion $x\in\quat$ is written as
\begin{align} \label{quaternion}
x = x_0+x_1 e_1+x_2 e_5 + x_3 e_6 
\end{align}
where $x_0,\,x_1,\,x_2,\, x_3 \in \mathbb{R}$, and we
have the usual rules of quaternionic multiplication
$e_1 e_5 = -e_5 e_1 = e_6$, etc. for the imaginary quaternionic units
$e_1, e_5$ and $e_6$.\footnote{The imaginary units are conventionally 
  designated as $i,j,k$ but here we prefer to use $e_1 \equiv i\,,\, 
  e_5 \equiv j$ and $e_6\equiv k$ because we want to reserve the 
  letter $i$ for the new imaginary unit used in the Cayley-Dickson doubling. 
  For the doubling $\quat = \cx\oplus i\cx$ we thus identify
  $i\equiv e_5$ and $x = (x_0+x_1 e_1)+ i (x_2-  x_3 e_1)$. Our notation 
  conforms with the one used in section~\ref{octsec} for the octonions, such
  that the obvious ${\rm D}_4$ roots inside ${\rm E}_8$ are simply  obtained from (\ref{D4})
  by multiplication with $e_7$.}
The norm is $|x|^2 \equiv x\bar{x} = x_0^2 + x_1^2 + x_2^2 + x_3^2$.

We are here interested in quaternionic integers. The {\em Hurwitz integers} 
(or just {\em Hurwitz numbers}) $\Hh$ are those quaternions 
(\ref{quaternion}) for which the coefficients 
$ x_0,\,
x_1,\,
x_2,\,
x_3$
are either all integers or all half-integers. They constitute 
a lattice $\Hh\subset\quat$ and form a non-commutative
ring of integers (actually, a `maximal order') inside $\quat$. 
The key feature is that $\Hh$ at the same time 
can be identified with the root lattice of the algebra ${\rm D}_4\equiv\mf{so}(4,4)$. Indeed,
the simple roots of ${\rm D}_4$, labeled according to the Dynkin diagram 
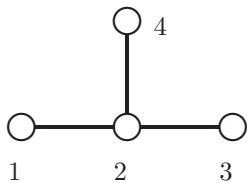
\begin{center}
\scalebox{1}{
\begin{picture}(348,60)
\put(125,0){$1$}
\put(165,0){$2$}
\put(205,0){$3$}
\put(180,55){$4$}
\thicklines
\multiput(130,20)(40,0){3}{\circle{10}}
\multiput(135,20)(40,0){2}{\line(1,0){30}}
\put(170,60){\circle{10}} 
\put(170,25){\line(0,1){30}}
\end{picture}}\end{center}
\vspace*{0.4cm}
can be identified with the following set of {\em Hurwitz units}:
\begin{align}\label{D4}
\ve_1 &= e_1, & 
\ve_2 &= \frac12(1-e_1-e_5- e_6), &
\ve_3&= e_5, &  \ve_4 &= e_6.
\end{align}
The ${\rm D}_4$ root lattice is spanned by integer linear combinations 
of these simple roots, 
and with the above choice of basis, the highest root 
$\theta = \ve_1 + 2\ve_2 + \ve_3 + \ve_4$ is indeed equal to $1$. 
The new feature is that such combinations also close under 
multiplication, thereby endowing the ${\rm D}_4$ root lattice with 
the structure of a non-commutative ring. The 24 roots of ${\rm D}_4$ are 
then identified with the 24 {\em units} in $\Hh$. Comparison 
with (\ref{sz}) shows that we now have the extra relation 
$s_1 s_3 s_4 = 1$ for the generating set (\ref{newevengenerators}), 
see (\ref{srelations}).

The elements in $W^+({\rm D}_4)$ are the words of even length where 
the letters are the fundamental Weyl reflections.
For example, those of length two act as
\begin{align}
x \mapsto \varepsilon_i (\overline{ \varepsilon_j \bar x \varepsilon_j }) \varepsilon_i=
\varepsilon_i (\bar{ \varepsilon}_j  x \bar{\varepsilon}_j) \varepsilon_i=
(\varepsilon_i \bar{ \varepsilon}_j)  x (\bar{\varepsilon}_j \varepsilon_i).
\end{align}
It follows \cite{Conway,Feingold:2008ih} that any element 
in $W^+({\rm D}_4)$ has the form of a combined left and right multiplication,
\begin{align} 
x \mapsto axb,
\end{align}
where $a$ and $b$ are (in general different) Hurwitz units, subject to
the constraint that the product $ab$ is an element of the quaternionic
group
\be\label{Q}
ab \in \cQ \equiv \{ \pm 1 ,\, \pm e_1,\,\pm e_5,\, \pm e_6 \}.
\ee
The restriction to the quaternionic group arises from triality of the ${\rm D}_4$ root system as the outer triality automorphisms are not elements of the Weyl group. In \cite{Feingold:2008ih} this result was generalized to the hyperbolic 
extension $W({\rm D}_4^{++})$ in a way that we will describe next.

\subsection{The even hyperbolic Weyl group $W^+({\rm D}_4^{++})$}
\label{sec:evenweyl}

Consider transformations acting on $\capX\in H_2(\quat)$ by
\begin{align}\label{evenquat}
\capX \; \rightarrow \; s(\capX) = S  \capX S^\dagger
\end{align}
with
\begin{align}  
S = \left( \begin{array}{cc} a & b \\
                             c & d \end{array}\right)
\qquad a,b,c,d\in \quat\,.
\end{align}
(Of course, these formulas trivially specialize to $\Div=\reals,\cx$.)
Such transformations preserve the norm of $\capX$ if and only if 
\begin{align}\label{detS}
\det (S S^\dagger) = 1\,.
\end{align}
This is true for all the normed division algebras up to and including the 
quaternions.\footnote{For real and complex numbers the determinant can 
 be factorized to give $|\det(S)|=1$. Making use of commutativity 
 then allows the reduction to $\det(S)=1$ in the projective version.} 
The associated continuous group of matrices will be denoted by 
\begin{align}\label{SL2H}
{\rm SL}(2,\quat)=\left\{S = \left( \begin{array}{cc} a & b \\
c & d \end{array}\right)\,:\,a,b,c,d\in \quat\,,\quad \det(SS^\dagger)=1\right\}.
\end{align}
Its projective version 
\begin{align}
{\rm PSL}(2,\quat):=
{\rm SL}(2,\quat)/\{ {\mathbf{1}}, - {\mathbf{1}}\big\}
\end{align}
is isomorphic to ${\rm SO}_0(1,5)$, see \cite{Su84,Baez}. Let us note the explicit 
expressions
\begin{align}\label{detSS}
\det (SS^\dagger) &=  |a|^2|d|^2+|b|^2|c|^2 - 
 2\,{\rm Re} \big(a\Bc d \Bb\big)    \nn\\
 &= |ad-bc|^2 - 2 \, {\rm Re} \big( a[\Bc,d] \Bb\big) 
\end{align}
and the explicit form of the (left and right) inverse matrix for $S\in SL(2,\quat)$
\begin{align}\label{Sinverse}  
S^{-1} = 
\left( \begin{array}{cc} |d|^2\Ba - \Bc d \Bb \, & \, |b|^2\Bc - \Ba b\Bd \\
         |c|^2 \Bb - \Bd c\Ba \, & \, |a|^2\Bd - \Bb a \Bc \end{array}\right)\,.
\end{align}
The discrete groups ${\rm SL}(2,\Hh)\subset {\rm SL}(2,\quat)$ and 
${\rm PSL}(2,\Hh)\subset {\rm PSL}(2,\quat)$ are obtained from the above 
groups by restricting the entries $a,b,c,d$ to be Hurwitz integers. 
From (\ref{detSS}) it then follows that in order for $S$ to be an 
element of ${\rm SL}(2,\Hh)$ or ${\rm PSL}(2,\Hh)$, the pairs $(a,b)$ and $(c,d)$ 
must each be {\em left coprime}, and the pairs $(a,c)$ and $(b,d)$ must 
each be {\em right coprime} (that is, share no common left or right 
factor $g\in\Hh$ with $|g|>1$). In fact, these conditions are not only necessary, but also sufficient (for some choice of the remaining two entries),
as we show in section \ref{we9we10-section}.

Straightforward computation using (\ref{capXform}) and (\ref{evenquat})
gives the transformations of the various components of $\capX$ 
\begin{align}\label{trmX}
s(x^+) &= \Ba a\, x^+ + \Bb b\, x^- + a x \Bb + b \bar x \Ba \,,   \nn\\
s(x)      &= a\Bc \, x^+ + b \Bd \, x^-  + a x \Bd  + b \bar x \Bc \,,\nn\\
s(x^-) &= \Bc c\, x^+ +  \Bd d\, x^- +  c x \Bd + d \bar x \Bc  \,.
\end{align}
This formula is, of course, also valid for $\Div=\reals,\cx$.
Likewise, it holds for general quaternionic matrices $S$ as long as 
these are invertible. 

The Weyl groups of the hyperbolic over-extensions of rank four algebras 
are discrete subgroups of ${\rm PSL}(2,\quat)$ (see \cite{Feingold:2008ih}
for details). 
For the 
hyperbolic over-extension ${\rm D}_4^{++}$ of ${\rm D}_4$ one finds
\begin{align}
W^+({\rm D}_4^{++}) \cong {\rm PSL}^{(0)}(2,\Hh) :=
\big\{ S\in {\rm PSL}(2,\Hh) \,|\, ad -bc \equiv 1 \,{\rm mod} \, \Cc \Big\}
\end{align}
where $\Cc$ is the (two-sided) commutator 
ideal $\Cc=\Hh[\Hh,\Hh]\Hh$
introduced in~\cite{Feingold:2008ih}.
In other words, the Weyl group of ${\rm D}_4^{++}$ is not simply ${\rm PSL}(2,\Hh)$:
the elements $S$ in ${\rm PSL}^{(0)}(2,\Hh)$ must satisfy the 
additional constraint $ad-bc\equiv 1$ modulo $\Cc$, which is an index 4 sublattice in $\Hh$. As shown 
in~\cite{Feingold:2008ih}, ${\rm PSL}^{(0)}(2,\Hh)$ is an index 3 subgroup of ${\rm PSL}(2,\Hh)$: if one extends 
the Weyl group by outer automorphisms related to the diagram automorphisms 
(`triality') one obtains all of ${\rm PSL}(2,\Hh)$, corresponding to all even 
symmetries of the ${\rm D}_4^{++}$ root lattice written in quaternionic coordinates.

Below, we will encounter two important subgroups
of ${\rm PSL}^{(0)}(2,\Hh)$. One is the (even) affine subgroup 
$W^+_{\text{aff}}\subset W^+({\rm D}_4^{++})$, which is the  
semi-direct product $W_{\text{fin}}\ltimes\cT$ of the finite Weyl group 
$W_{\text{fin}}=W^+({\rm D}_4)$ and the abelian group $\cT$ of affine 
translations. This is a maximal parabolic subgroup of $W^+({\rm D}_4^{++})$
whose elements in the quaternionic matrix representation are given by
\begin{align}\label{Par1}
S =
\begin{pmatrix}
a & u\\
0 & b
\end{pmatrix}\;\; ,  
\end{align}
where $a,b$ are unit quaternions with $ab\in\cQ$ and $u\in\Hh$ is an 
arbitrary Hurwitz integer. The other important subgroup is the
translation subgroup $\cT$ itself, consisting of the matrices
\begin{align}\label{Par2}
T_u =
\begin{pmatrix}
1 & u\\
0 & 1
\end{pmatrix}\;\; ,  
\end{align}
where again $u\in\Hh$ is an arbitrary Hurwitz number. The following
little lemma will be useful in section~\ref{sec:eisen}.

\begin{lemma}\label{hurwlemma}
The left coprime pairs $(c,d)\in \Hh^2$
parametrize the coset spaces {\em $W^+({\rm D}_4^+)\backslash W^+({\rm D}_4^{++})\equiv 
W^+_{\text{aff}}\backslash W^+_{\text{hyp}}$}. The equivalence classes $(c,d)\sim (ec,ed)$ for a unit $e\in\Hh$ uniquely parametrize the cosets.
\end{lemma}

\Pf Suppose that the left coprime lower entries $(c,d)$ are given 
and that  $a,b$ and $\tilde{a},\tilde{b}$ are two different pairs of 
Hurwitz numbers completing $c,d$ to two different matrices $S$ 
and $\tilde{S}$ in ${\rm PSL}^{(0)}(2,\Hh)$. Then an easy computation using
(\ref{Sinverse}) shows that
\be
\tilde{S}S^{-1} = \begin{pmatrix}
                         q & *\\
                         0 & 1
                  \end{pmatrix}\,,
\ee
where $q\in \cQ$ is in the quaternionic group (\ref{Q}).
Hence all such matrices are related by an upper triangular element 
of ${\rm PSL}^{(0)}(2,\Hh)$ of type (\ref{Par1}), that is, an element of $W^+_{\text{aff}}$. Since (\ref{Par1}) also allows for left multiplication of the pair $(c,d)$ by a Hurwitz unit we arrive at the claim.
\qed

\subsection{Action of ${\rm PSL}^{(0)}(2,\Hh)$ on $\cH(\quat)$}

We now wish to interpret the even Weyl group 
$W^+({\rm D}_4^{++})\equiv {\rm PSL}^{(0)}(2,\Hh)$ as a 
modular group $\Gamma$ acting on the quaternionic upper half plane
$\cH(\quat)$. To exhibit the nonlinear `modular' action of $\Gamma$ 
we map the forward unit hyperboloid $||\capX||^2=-1$ in $\reals^{1,5}$ 
to the upper half plane $\cH(\quat)$ by means of the projection 
(\ref{uhpco}). Accordingly, we consider (\ref{z}), but now specialize 
to $u\in\quat$. Using the formulas of the previous section, 
in particular  (\ref{uhpco}) and (\ref{trmX}), we obtain 
\begin{align}\label{uvtrm}
v' = \frac{v}{D} \;\;,\quad
u' = \frac1{D} \Big[  (au+b)(\Bu\Bc+\Bd) + a\Bc\, v^2 \Big]
\end{align}
with $z' = u' + iv'\equiv s(z) \equiv s(u+iv)$, where
\begin{align}
D\equiv D(S,u,v) := |cu + d|^2 + |c|^2 v^2
= |cz + d|^2
\end{align}
and (\ref{CD1}) has been used. Observe that
now $cz+d \in \oct$. The transformations (\ref{uvtrm})
can be combined into a single formula
\begin{align}\label{z'}
z'  =
\frac{ (au+b)(\Bu\Bc+\Bd) + a\Bc \, v^2  + iv}{|cz + d|^2}\;\;\;, \qquad
\begin{pmatrix}
a & b\\
c & d
\end{pmatrix}\in {\rm PSL}^{(0)}(2,\Hh),
\end{align}
which is our expression for the most general quaternionic
modular transformation. Here $z'$ is manifestly in the upper half plane
$\cH(\quat)$, and  also reduces to the standard formula $z' = (az+b)(cz+d)^{-1}$ for 
the commutative cases $u\in\reals$ (or $\cx$) and $a,b,c,d\in\reals$ 
(or $\cx$).
However, (\ref{z'}) fails for $\Div =\oct$ because of 
non-associativity. One can check explicitly that the transformations 
(\ref{uvtrm}) and the identifications (\ref{uhpco}) are consistent 
with (\ref{trmX}). 

As a special case we have the modular transformations (\ref{Par1})
corresponding to the even affine Weyl group.
The induced modular action of the most general element
$\gamma\in\Gamma_\infty \equiv W_{\text{aff}}^+\equiv 
W_{\text{fin}}^{+}\ltimes \Hh$ is 
\begin{align}\label{affine}
\gamma  \; :\;
z \rightarrow \gamma(z) = aub  + y\bar{b} + iv
\end{align}
with units $a,b$ such that $ab\in\cQ$ (cf.~(\ref{Q}))
and $y\in\Hh$ an arbitrary Hurwitz integer. This transformation 
corresponds to a finite (even) Weyl transformation followed by 
a constant shift along the root lattice $\Hh$. In particular, 
it follows that such transformations leave invariant the `cusp' $z=i\infty$ 
in $\cH(\quat)$ --- in analogy with the action of the usual shift 
matrix $T$ for the real modular group ${\rm PSL}(2,\ints)$ (in that case, 
$W_{\text{fin}}^+$ is trivial). As we explained, this is in accord 
with the fact that the affine Weyl group $W_{\text{aff}}$ is the 
subgroup of $W_{\text{hyp}}$ leaving invariant the affine null root 
$\delta$ \cite{GN}.

\section{${\rm E}_8$ and its hyperbolic extension}
\label{octsec}

We now turn to the largest exceptional algebra ${\rm E}_8$ and its affine and hyperbolic extensions ${\rm E}_9\equiv {\rm E}_8^+$ and ${\rm E}_{10}\equiv {\rm E}_8^{++}$. When extending the results of the previous section to $\Div=\oct$ a 
main obstacle is the non-associativity of the octonions. We have the abstract
isomorphism~\cite{Feingold:2008ih}
\begin{align}
W^+({\rm E}_{10}) \cong {\rm PSL}(2,\Oo)
\end{align}
and we will exhibit a presentation of this group in terms of cosets with respect to the cusp stabilizing affine Weyl group $W^+(E_9)$; a presentation that will prove particularly useful when constructing non-holomorphic automorphic forms in section~\ref{sec:quatauto}. We emphasize that realizing part of this group
in terms of $2 \times 2$ matrices with integer octonionic (`octavian') 
entries is certainly not sufficient, because such matrices violate
associativity and thus cannot by themselves define any group.
Nevertheless, the generating transformations (\ref{sz}) for 
$\Gamma\equiv {\rm PSL}(2,\Oo)$ are still valid in the form given there:
by the alternativity of the octonions no parentheses need to be
specified for nested products involving only two different octonions. 
On our way to describing $W^+({\rm E}_{10})$ we will also obtain formula (\ref{newwe8}) below  that gives a rather explicit expression
for the action of an ${\rm E}_8$ Weyl group element in terms of unit octavians and
the ${\rm G}_2(2)$ automorphisms of the octavians.

\subsection{${\rm E}_8$ and the integer octonions}

As a basis for the octonion algebra one usually takes the real number 
$1$ and seven `imaginary' units $e_1,\,e_2,\,\ldots,\,e_7$ that 
anticommute and square to $-1$. One then chooses three of them 
to be in a quaternionic subalgebra. For instance, taking $(e_1,\,e_5,\,e_6)$
as a basic quaternionic triplet as in section~\ref{sec:d4sec}, the octonion multiplication 
has the following properties:
\begin{align}\label{octmult}
e_ie_j=e_k \quad \Rightarrow \quad e_{i+1}e_{j+1}=e_{k+1},\quad
e_{2i}e_{2j}=e_{2k},
\end{align} 
where the indices are counted mod 7. Using these rules, the whole 
multiplication table can be derived from only one equation, for 
example $e_1e_5 = e_6$.

The root lattice of ${\rm E}_8$ can be identified with the non-commutative 
and non-associative ring of {\em octavian integers}, or simply {\em octavians} $\Oo$ (for their construction, see
\cite{Conway}). Within
this ring, the 240 roots of ${\rm E}_8$ are then identified with the 
set of {\em unit octavians}, which are the 240 invertible elements in 
$\Oo$. More specifically, 
the octavian units can be divided into
\begin{align}
2\ &\textit{real numbers:} &&\pm1,\nn\\
112\ &\textit{Brandt numbers:} &&\tfrac12(\pm1\pm e_i\pm e_j\pm e_k), \nn\\
126\ &\textit{imaginary numbers:} &&\tfrac12(\pm e_m\pm e_n\pm e_p\pm e_q),\quad \pm e_r,
\end{align}
where
\begin{align}
ijk\ =&&124&,&137&,&156&,&236&,&257&,&345&,&467,&\nn\\
mnpq\ =&&3567&,&2456&,&2347&,&1457&,&1346&,&1267&,&1235,&
\end{align}
and $r=1,2,\ldots,7$.~\footnote{Imaginary units and Brandt numbers are
   called "eyes" and "arms" in \cite{Conway}.}  
The lattice $\Oo$ of octavians is the 
set of all integral linear combinations of these unit octavians.

Recall that the automorphism group ${\text{Aut}}\,\oct$ of the octonionic multiplication 
table (\ref{octmult}) is the exceptional group ${\rm G}_2$ \cite{Baez}. For
the ring $\Oo$ we have the analogous discrete result \cite{Conway}
\begin{align}
{\text{Aut}}\, \Oo &= {\rm G}_2(2), & |{\rm G}_2(2)| &= 12\, 096,
\end{align}
where ${\rm G}_2(2)$ is the finite (but not simple) group ${\rm G}_2$ over the 
field ${\mathbb{F}}_2$ (a Chevalley group).

We choose the simple roots of ${\rm E}_8$, according to the Dynkin diagram
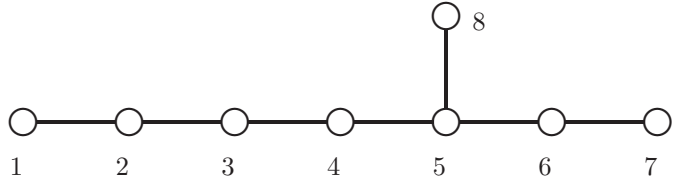
\begin{center}
\scalebox{1}{
\begin{picture}(260,60)
\put(5,-10){$1$}
\put(45,-10){$2$}
\put(85,-10){$3$}
\put(125,-10){$4$}
\put(165,-10){$5$}
\put(205,-10){$6$}
\put(245,-10){$7$}
\put(180,45){$8$}
\thicklines
\multiput(10,10)(40,0){7}{\circle{10}}
\multiput(15,10)(40,0){6}{\line(1,0){30}}
\put(170,50){\circle{10}} \put(170,15){\line(0,1){30}}
\end{picture}}\end{center}
\vspace*{0.4cm}
as the following unit octavians, 
\begin{align}
\varepsilon_1 &= \tfrac12 (1-e_1-e_5-e_6),&
\varepsilon_2 &= e_1,\nn\\
\varepsilon_3 &= \tfrac12 (-e_1-e_2+e_6+e_7),&
\varepsilon_4 &= e_2,\nn\\
\varepsilon_5 &= \tfrac12 (-e_2-e_3-e_4-e_7),&
\varepsilon_6 &= e_3,\nn\\
\varepsilon_7 &= \tfrac12 (-e_3+e_5-e_6+e_7),&
\varepsilon_8 &= e_4.\label{e8roots}
\end{align}
\noindent
The identification (\ref{e8roots}) 
of the simple roots of ${\rm E}_8$ as unit octavians is not unique.
The main advantage of the choice of basis here (similar to~\cite{KocaKarsch89} but different from the one
in \cite{Feingold:2008ih}) is that the 126 roots of the ${\rm E}_7$ subalgebra 
are given by {\em imaginary} octavians (defined to obey $e^2=-1$), while 
the 112 roots corresponding to the two $\bf{56}$ representations of ${\rm E}_7$ 
in the decomposition of the adjoint $\bf{248}$ of ${\rm E}_8$ are given 
by {\em Brandt numbers}, and the highest root is $\theta=1$ in accord 
with our general conventions. By definition, a Brandt number $a$ 
obeys $a^3=\pm 1$. As explained in \cite{Zo35,Conway}, for such numbers (and
only for them) the map $x\rightarrow axa^{-1}$  is an automorphism of 
the octonions, see also appendix \ref{autoctapp}. 

Note that only the leftmost node of the Dynkin diagram corresponds 
to a Brandt number, while all others correspond to imaginary unit octavians. 
Also, 
the root basis for the obvious ${\rm D}_4$ subalgebra inside ${\rm E}_8$ is simply obtained 
from the one of section~\ref{quat-section} by multiplying 
$\varepsilon_4, \varepsilon_5, \varepsilon_6$
and $\varepsilon_8$ with $e_7$.
 
Using the scalar product (\ref{ab}) it is straightforward to reproduce 
from (\ref{e8roots}) the ${\rm E}_8$ Cartan matrix via
\begin{align}
A_{ij} = 2(\varepsilon_i , \varepsilon_j).
\end{align}

We saw in section~\ref{quat-section} that any element in the finite Weyl group 
$W^+({\rm D}_4)$ can be written as a bi-multiplication
\begin{align}
x \mapsto axb, \label{d4transformation}
\end{align}
where $a$ and $b$ are Hurwitz units. But the number of different such 
transformations (= 12 $\times$ 24) exceeds the order of $W^+({\rm D}_4)$
(which is = 96). This is why, in order to obtain only elements in 
$W^+({\rm D}_4)$, we needed to impose the extra condition 
$ab\in \cQ$. For $W^+({\rm E}_8)$ we would 
like to generalize the expression (\ref{d4transformation}) to unit 
{\it octavians} $a$ and $b$. Because of the non-associativity, we 
have to place parentheses, for instance by taking
\begin{align}
x \mapsto (ax)b. \label{e8transformation}
\end{align}
Since the Dynkin diagram of ${\rm E}_8$ does not have any symmetries 
(unlike the one of ${\rm D}_4$) all such transformations do belong to 
$W^+({\rm E}_8)$. But counting these expressions, we find that their 
number (modulo sign changes $(a,\,b) \leftrightarrow (-a,\,-b)$) 
now is {\it smaller} than the order of the even Weyl group, which is
\begin{align}
|W^+({\rm E}_8)|=120 \times 240 \times 12\,096. \label{e8evenweylorder}
\end{align}
The extra factor $12\,096$ is exactly the order of the automorphism 
group of the octavians. This observation, made in \cite{Feingold:2008ih} 
(see also \cite{Ramond}), suggests that the expression 
(\ref{e8transformation}) should be generalized to
\begin{align}
x \mapsto (a\varphi(x))b, \label{e8transformationauto}
\end{align}
where $\varphi$ is an arbitrary automorphism of the octavians. 
However, it was noted in \cite{Feingold:2008ih} that the number of
{\it independent} such expressions is in fact smaller than 
$120 \times 240 \times 12\,096$,
since the transformation (\ref{e8transformation}) is already an automorphism 
if $a$ is a Brandt number (that is, $a^3= \pm 1$) and $b=\pm \bar a$ (see above and appendix \ref{autoctapp}).
In section~\ref{e7e8} we will explain how the formula (\ref{e8transformationauto}) can be modified in these cases, so that it indeed expresses all even Weyl transformations of ${\rm E}_8$.
As a first step we consider the ${\rm E}_7$ subalgebra obtained by deleting the leftmost node in the Dynkin diagram and thus corresponding to imaginary octavians.

\subsection{The even Weyl group of ${\rm E}_7$} \label{evene7weylgroup}
In this subsection we use $a_\sB$ to denote bimultiplication from left and right with the same octavian $a$:
\begin{align}
a_\sB :\quad \Oo \to \Oo, \quad x \mapsto a x a.
\end{align}
We will show that
$W^+({\rm E}_7)$ is the set of all transformations 
\begin{align} \label{e7weylelement}
g_{\sB} h_\sB \varphi: \quad  \Oo \to \Oo,\quad
x \mapsto g(h \,\varphi(x) h) g,
\end{align}
where $g$ and $h$ are imaginary unit octavians and $\varphi$ is an automorphism of the octavians.
It follows from (\ref{wx}), and the identification
(\ref{e8roots}) of the simple ${\rm E}_7$ roots as imaginary unit octavians, that
the generators of $W^+({\rm E}_7)$ have this form.
To show that this holds for all elements in $W^+({\rm E}_7)$
we thus have to show that the composition of two such transformations 
$g_{1\sB} h_{1\sB} \varphi_1$ and ${g}_{2\sB} {h}_{2\sB} \varphi_2$
again can be written in the form (\ref{e7weylelement}).
From the definition of an automorphism we have
\begin{align} \label{product2}
(g_{1\sB} h_{1\sB} \varphi_1) (g_{2\sB} h_{2\sB} \varphi_2)
&= g_{1\sB} h_{1\sB}\, \varphi_1(g_2)_\sB\, \varphi_1(h_2)_\sB \, \varphi_1 \varphi_2\nn\\
&= (e_\sB f_\sB f_\sB e_\sB) g_{1\sB} h_{1\sB}\, \varphi_1(g_2)_\sB\, \varphi_1(h_2)_\sB \varphi_1 \varphi_2,
\end{align}
where we in the last step have inserted the identity map in form of $e_\sB f_\sB f_\sB e_\sB$ for two arbitrary imaginary unit octavians $e$ and $f$.
Any unit octavian can be written as a product of two imaginary unit octavians, so we can choose
$e$ and $f$ in (\ref{product2}) such that
\begin{align}
fe = \pm (  ( \,\varphi_1(h_2) \,\varphi_1(g_2) ) h_1 ) g_1,
\end{align}
or equivalently,
\begin{align}
((((fe) g_1) h_1 ) \,\varphi_1(g_2)) \,\varphi_1(h_2)=\pm1.
\end{align}
Then, according to Corollary \ref{thecorollary}, 
the map
\begin{align}
\varphi_3 \equiv f_\sB e_\sB g_{1\sB} h_{1\sB} \,\varphi_1(g_2)_\sB \,\varphi_1(h_2)_\sB
\end{align}
is an automorphism.
Inserting this into (\ref{product2}) we obtain
\begin{align} \label{product3}
(g_{1\sB} h_{2\sB} \varphi_1) (g_{2\sB} h_{2\sB} \varphi_2)
&= (e_\sB f_\sB f_\sB e_\sB) g_{1\sB} h_{1\sB} \,\varphi_1(g_2)_\sB \,\varphi_1(h_2)_\sB \varphi_1 \varphi_2\nn\\
&= e_\sB f_\sB ( f_\sB e_{\sB} g_{1\sB} h_{1\sB} \,\varphi_1(g_2)_\sB \,\varphi_1(h_2)_\sB ) \varphi_1 \varphi_2\nn\\
&= e_\sB f_\sB (\varphi_3 \varphi_2 \varphi_1).
\end{align}
Since $\varphi_3,\,\varphi_1$ and $\varphi_2$ all are automorphisms, their product $\varphi_3 \varphi_2 \varphi_1$ is an automorphism as well.

Conversely, any 
transformation (\ref{e7weylelement}) is an 
isometry of the root system.
Since the group of isometries of the root system of a finite algebra is the 
semidirect product of the Weyl group and the symmetry group of the 
Dynkin diagram (which is trivial for ${\rm E}_7$), it follows that any transformation
(\ref{e7weylelement}) belongs to the Weyl group. Furthermore, it belongs to the {\em even} Weyl group since it does not involve the conjugate $\bar x$.
 
We have thus proven that 
$W^+({\rm E}_7)$ is the set of 
all expressions (\ref{e7weylelement}). But there is some redundancy in this set.
If $g_1 h_1 = g_2 h_2$, then the transformation
\begin{align}
\varphi_3=h_2 g_2 g_1 h_1
\end{align}
is an automorphism according to Corollary \ref{thecorollary}, and we have 
\begin{align}
g_1 h_1 \varphi_1 = g_2 h_2 \varphi_2
\end{align}
where $\varphi_2=\varphi_3\varphi_1$.
Thus only the product $\pm gh$ matters modulo automorphisms $\varphi$, and this 
product can be any unit octavian.  Put differently, any element of $W^+({\rm E}_7)$
is characterized by an automorphism in $G_2(2)$ and a unit octavian $a$ which
must then be factorized as $gh$ with imaginary units $g$ and $h$, unless $a$
is itself imaginary. Accordingly, the order of $W^+({\rm E}_7)$ is 
$120\times 12\,096$.

\subsection{From $W({\rm E}_7)$ to $W({\rm E}_8)$} \label{e7e8}

Having understood the Weyl group of ${\rm E}_7$ the step to ${\rm E}_8$ is more conceptual than technical -- it is an application of the orbit-stabilizer theorem, which is valid for any group $\mathcal{G}$ acting on a set $\mathcal{X}$.
If $\mathcal{H} \subset \mathcal{G}$ is the subgroup stabilizing an element $x \in \mathcal{X}$, the theorem says that, for any 
$g \in \mathcal{G}$, the coset $g\mathcal{H}$ is the set of all elements in $\mathcal{G}$ that take $x$ to $gx$.

The Weyl group of ${\rm E}_7$ is the subgroup of 
$W({\rm E}_8)$ that stabilizes the highest root of ${\rm E}_8$. 
Applying the orbit-stabilizer theorem, it is enough to find one transformation $w$ in $W({\rm E}_8)$ for each element $x$ in the orbit $W({\rm E}_8)\,\theta$,
such that $w(\theta)=x$. Then the set of {\em all} such transformations 
is the coset $w\,W({\rm E}_7)$.
The root system of ${\rm E}_8$ is one single orbit under $W({\rm E}_8)$ and the highest root $\theta$ is 1 in our identification with the octavians. Thus we need a transformation in $W({\rm E}_8)$ for each unit octavian $b$, such that $\theta=1$ is mapped to $b$. But this is easy to find, we just take the transformation to be right (or left) multiplication by $b$. It follows if we linearize the composition property (\ref{compprop}), or alternatively by (\ref{ab}) together with the first of the Moufang identities (\ref{moufang}), that
\begin{align}
(xb,\,yb)=|b|^2(x,\,y),
\end{align}
so right multiplication by $b$ is an isometry if $b$ has unit length. Like for ${\rm E}_7$ any isometry of the root system is an element of the Weyl group, and right multiplications of unit octavians belong to the even Weyl group $W^+({\rm E}_8)$ since they do not involve conjugation.

Thus any element of $W^+({\rm E}_8)$ can be written as
\begin{align} \label{e8weylelement}
x \mapsto (f(e \,\varphi (x)\,e)f)b,
\end{align}
where $e$ and $f$ are imaginary unit octavians, $b$ an arbitrary unit octavian, and $\varphi$ an automorphism of the octavians. Using the third of the Moufang identities (\ref{moufang}) we can write this as
\begin{align}\label{newwe8}
x &\mapsto (f(e\,\varphi (x)\,e)f)b\nn\\
&\quad\,=f((e\,\varphi(x)\,e)(fb))\nn\\
&\quad\,=f(e\,(\varphi(x)\,(e(fb))))\nn\\
&\quad\,=f(e\,(\varphi(x)\,d))
\end{align}
where $d=e(fb)$. This differs from the formula (\ref{e8transformationauto})
suggested in \cite{Feingold:2008ih}, only by the factorization of $a$ into two imaginary unit octavians $e$ and $f$. 
In fact, the factors $e$ and $f$ do not need to be imaginary but can also be {\em real}, that is $\pm1$. Then if $a$ in  
(\ref{e8transformationauto}) is imaginary, the factorization is trivial; we can just take $e=a$ and $f=1$, and 
(\ref{newwe8}) reduces to (\ref{e8transformationauto}). 
Only when $a$ is a Brandt number we need to replace (\ref{e8transformationauto}) by (\ref{newwe8}), where we can choose any factorization of $a$ into two imaginary unit octavians $e$ and $f$.

When we consider ${\rm E}_8$ as a subalgebra of
${\rm E}_{10}={\rm E}_8^{++}$, the Weyl group $W({\rm E}_8)$ acts in a non-trivial way on the simple root $\alpha_0$ (and trivially on $\alpha_{-1}$). We can thus extend the action of $W({\rm E}_8)$ to the whole root space of 
${\rm E}_{10}$. The Weyl transformation (\ref{newwe8}) can then be written as\footnote{The fact that the 
formula `$\capX\mapsto S\capX S^\dagger$' requires refinement 
by placing parentheses inside the matrix elements was 
anticipated in~\cite{Feingold:2008ih}.}
\begin{align}
\begin{pmatrix}
x^+ & x \\
\bar x &x^-
\end{pmatrix}\mapsto
\begin{pmatrix}
x^+ & f(e\,(\varphi(x)\,d))\\
((\bar d\,\varphi(\bar x))\bar e)\bar f &x^-
\end{pmatrix}
\end{align}
or, in more compact form,
\begin{align}
\capX\mapsto
F(E(D^\dagger\varphi(\capX)D)E^\dagger)F^\dagger, \label{we8matrix}
\end{align}
where
\begin{align} \label{chidef}
D&=
\begin{pmatrix}
1 & 0 \\
0 & d
\end{pmatrix},&
E&=
\begin{pmatrix}
e & 0 \\
0 & 1
\end{pmatrix},&
F&=
\begin{pmatrix}
f & 0 \\
0 & 1
\end{pmatrix}.
\end{align}

\subsection{From $W({\rm E}_8)$ to $W({\rm E}_9)$}

As mentioned in section~\ref{sec:affweyl}  (see also \cite{Ka90}), to any point $x$ on the root lattice 
of the finite Kac-Moody algebra ${\rm E}_8$ we can associate a {\it translation} 
$t_x$ that acts on the root space of ${\rm E}_{10}$ as 
\begin{align}
t_y : \capX \mapsto \capX + 2(\capX,\delta)y + 
  2(\capX,\,y)\delta-2|y|^2(\capX,\,\delta)\delta.
\end{align}
where $\delta$ is the null root (\ref{nullroot}). Then $t_x t_y=t_{x+y}$ 
for any $x$ and $y$, and all translations form a free group of translations
$\mathcal{T}\equiv \Oo $ of rank eight. This group is a normal subgroup 
of the Weyl group of ${\rm E}_9$, which can then be written as a semidirect 
product
\begin{align}
W({\rm E}_9)= W({\rm E}_8) \ltimes \mathcal{T} \equiv W({\rm E}_8) \ltimes \Oo.
\end{align}
Furthermore, the translations belong to the even Weyl group, so we have
\begin{align}
W^+({\rm E}_9)=W^+({\rm E}_8)  \ltimes \mathcal{T}.
\end{align}
When we identify the root space of ${\rm E}_{10}$ with the Jordan algebra 
$H_2(\oct)$, we can write the action of $t_y$ as
\begin{align} \label{transdef0}
t_y : \capX \mapsto  T_y \, \capX \, T_y^\dagger
\end{align}
where
\begin{align}
\label{transdef}
T_y = 
\begin{pmatrix}
1 & y\\
0 & 1
\end{pmatrix}.
\end{align}
For such matrices we indeed need not worry about non-associativity because
\begin{align}
({T_y} \capX ){T_y}^\dagger={T_y} (\capX {T_y}^\dagger)&=
\begin{pmatrix}
x^+ + y\bar x + x \bar y + yx^- \bar y & x+yx^-\\
\bar x+x^-\bar y & x^-
\end{pmatrix}\nn\\&=
\begin{pmatrix}
x^+ +2(x,\,y) +x^- |y|^2& x+x^-y\\
\bar x+x^-\bar y & x^-
\end{pmatrix}=t_y(\capX).
\end{align}
where $\delta$ is as in (\ref{nullroot}).
It follows that any element in $W^+({\rm E}_9)$ can be written
in the form (\ref{we8matrix}), but now with
\begin{align}
D&=
\begin{pmatrix}
1 & 0 \\
c & d
\end{pmatrix}
\end{align}
where $c$ is an arbitrary octavian,
so that $D^\dagger$ is an upper triangular matrix. 
Since at least one of the diagonal entries in $D$ is equal to 1 the 
expression (\ref{we8matrix}) is still well defined. 

\subsection{From $W({\rm E}_9)$ to $W({\rm E}_{10})$} \label{we9we10-section}

The tool for describing $W({\rm E}_{10})$ in terms of $W({\rm E}_9)$ is the orbit-stabilizer theorem, just like in going from $W({\rm E}_7)$ to $W({\rm E}_8)$. We use the fact that $W^+({\rm E}_9)$ stabilizes the null root $-\delta$ of (\ref{nullroot}). Therefore the left cosets $W^+({\rm E}_{10})/W^+({\rm E}_9)$ are in bijection with the orbit of $-\delta$ under $W^+(E_{10})$. 
In this section we will see that the orbit of $-\delta$ under $W^+(E_{10})$ in turn can be parametrized by pairs of right coprime octavians, if we define coprimality for octavians in an appropriate way based on the Euclidean algorithm. (This is not the only possible generalization of the notion of coprimality when going from the Hurwitz numbers to the non-associative octavians --- see below).
We first recall the validity of the (right) Euclidean algorithm for octavians (see e.g.~\cite{Conway} for a proof).

\begin{thm} 
Let $a$ and $c$ be two octavians. Then, for some integer $n\geq0$, there exist octavians 
$q_1$, $q_2$, $\ldots$, $q_n$, $q_{n+1}$ and 
$r_0$, $r_1$, $\ldots$, $r_n$,
such that
\begin{align}
a &= q_1 c - r_1,\nn\\
c &= q_2 r_1 - r_2\nn\\
r_1 &= q_3 r_2 - r_3,\nn\\
&\cdots\nn\\
r_{n-2} &= q_{n} r_{n-1} - r_{n},\nn\\
r_{n-1} &= q_{n+1} r_{n} \label{euclalg}
\end{align}
and $|r_1|>\cdots>|r_n|>0$.\footnote{In this way of writing the Euclidean algorithm, we have changed the sign of the remainders $r_i$ compared to other authors. This choice of sign proves to be more convenient in the analysis to follow.}
\end{thm}
We now define $a$ and $c$ to be right coprime if (for some choice of possible $n$ and 
$q_1,q_2 \ldots, q_{n+1}$), the last non-vanishing remainder is a unit, $|r_n|=1$.
For Hurwitz numbers this definition is equivalent to the one that we gave in
section \ref{sec:evenweyl}: Two Hurwitz numbers are left (right) coprime if and only if they share no common left (right)
factor $g\in\Hh$ with $|g|>1$.
For octavians the two definitions are no longer equivalent, due to non-associativity. One can construct counterexamples of $a$ and $c$ such that they are right coprime but still have a non-trivial common right divisor.\footnote{Consider for example $a=e_1+e_2$, $c=e_1+e_3$. Then one can write $a=q_1 c - r_1$ with
\begin{align}
q_1 = \frac12\left(2+e_1+e_4-e_5-e_7 \right)\,,\quad
r_1 = \frac12\left( -1+e_1+e_3+e_7\right),\nn
\end{align}
so they right co-prime in the sense of (\ref{euclalg}). At the same time, they have the common right divisor $g=1+e_1$.} 
We have not been able to decide on the converse statement, i.e., if having only trivial common right divisor implies being right coprime in the sense of the Euclidean algorithm. Left coprimality is defined analogously as we will spell out below.

With the above definition of coprimality of octavians we have
the following lemma.
\begin{lemma} \label{orbitlemma}
Let $\cO_{-\delta}$ be the orbit of $-\delta$ under the action of $W^+({\rm E}_{10})$. Then
\be\label{deltaorbit}
\cO_{-\delta} = \left\{ \begin{pmatrix}
|a|^2 & a\bar{c}\\
c\bar{a} & |c|^2 
\end{pmatrix}
\,:\, \text{\rm $a$ and $c$ are right coprime octavians}\right\}\,.
\ee
\end{lemma}

\Pf Let $\mathcal{S}$ be the right hand side of (\ref{deltaorbit}). By examining the action of the generators (\ref{evenweylgen}) of $W^+({\rm E}_{10})$ one finds that $\mathcal{S}$ is preserved under the action of $W^+({\rm E}_{10})$. Since $-\delta$ itself is contained in $\mathcal{S}$ we conclude that $\cO_{-\delta}$ is contained in $\mathcal{S}$. 

To prove the other inclusion, we need to find a $W^+({\rm E}_{10})$ element $w_{a,c}$ 
associated to any pair of right coprime octavians $a$ and $c$ such that $w_{a,c}(-\delta)$ 
equals the matrix on the right hand side of (\ref{deltaorbit}). 
Since $a$ and $c$ are right coprime, there exist octavians $q_1, q_2, \ldots, q_{n+1}$ 
such that the remainders $r_1,r_2,\ldots,r_n$, defined by (\ref{euclalg}), satisfy 
$|r_1|>\cdots>|r_n|>0$ and $|r_n|=1$.
We now consider the $W^+({\rm E}_{10})$ element
\be\label{wac}
w_{a,c} = t_{q_1}\circ s_{-1} \circ \dots \circ t_{q_{n+1}}\circ s_{-1}\circ u_{r_n}
\ee
where $t_q$ has been defined in (\ref{transdef0}), $s_{-1}$ is as in (\ref{evenweylgen}) and $u_{r_n}$ is given by
\be\label{rotgen}
u_{r_n}(X) = U_{r_n} X U_{r_n}^\dagger\quad\quad\text{with}\quad
U_{r_n} = \begin{pmatrix}
r_n & 0 \\
0 & \bar{r}_n
\end{pmatrix}.
\ee
Thus $u_{r_n}$ is an element of $W^+({\rm E}_{8})\subset W^+({\rm E}_{10})$ if and only if $|r_n|=1$. The notation $w_{a,c}$ is not completely accurate as the element also depends on the choice of Euclidean algorithm decomposition, but all (presently) important features of $w_{a,c}$ depend only on $a$ and $c$. We 
calculate $w_{a,c}(-\delta)$ by induction. One first verifies that
\be
(t_{q_{i+1}}\circ s_{-1}) \begin{pmatrix}
|r_i|^2 & r_i \bar{r}_{i+1}\\
r_{i+1} \bar{r}_i & |r_{i+1}|^2
\end{pmatrix}
=\begin{pmatrix}
|r_{i-1}|^2 & r_{i-1} \bar{r}_{i}\\
r_{i} \bar{r}_{i-1} & |r_{i}|^2
\end{pmatrix}
\ee
for all $i=0,\ldots,n$, if we let $r_{-1}$, $r_0$ and $r_{n+1}$ be equal to $a$, $c$ and $0$, respectively.
Since trivially 
\be\label{rotdelta}
u_{r_n} (-\delta) = \begin{pmatrix}
1 & 0\\
0 &0 \end{pmatrix}
=\begin{pmatrix}
|r_n|^2 & r_n \bar{r}_{n+1}\\
r_{n+1} \bar{r}_n & |r_{n+1}|^2
\end{pmatrix}
\ee
we conclude that
\be
w_{a,c} (-\delta) = \begin{pmatrix}
|r_{-1}|^2 & r_{-1}\bar{r}_0\\
r_0 \bar{r}_{-1} & |r_0|^2 \end{pmatrix}
=\begin{pmatrix}
|a|^2 & a\bar{c}\\
c \bar{a} & |c|^2\end{pmatrix}.
\ee
This shows that all elements of $\mathcal{S}$ are contained in $\cO_{-\delta}$, completing the proof of the lemma.
\qed

\noindent
By the orbit-stabilizer theorem 
this lemma
implies that
$W^+(E_{10})$ can be written as
\begin{align}
\label{WE10WE9}
W^+({\rm E}_{10}) = \bigcup_{a,c\in\Oo\,\, \text{\rm right coprime}} 
w_{a,c} 
W^+({\rm E}_9)
\end{align}

By repeating the same arguments leading to (\ref{WE10WE9}) but using a right action with the inverse elements we find that
$W^+(E_{10})$ is 
not only the union of left cosets (\ref{WE10WE9})
but also the union of the right cosets
\begin{align}\label{WE9WE10}
W^+({\rm E}_{10}) = \bigcup_{c,d\in\Oo\,\, \text{left coprime}}  
W^+({\rm E}_9)
\tilde{w}_{c,d}\,.
\end{align}
where we sum now over {\em left} coprime octavians and
\begin{align}\label{wcd}
\tilde{w}_{c,d} = u_{\bar{r}_n} \circ s_{-1}\circ t_{q_{n+1}}  \circ\dots\circ s_{-1}\circ t_{q_1}.
\end{align}
These octavians $q_1,\ldots q_{n+1}$ and $r_n$ appearing here are the elements in the left Euclidean algorithm
\begin{align}\label{leuclalg}
 d &= c q_1 - r_1 \,,\nn\\
 c &= r_1 q_2- r_2 \,,\nn\\
 r_1 &= r_2 q_3 - r_3\,,\nn\\
&\dots\nn\\
 r_{n-2} &= r_{n-1} q_n- r_n\,,\nn\\
 r_{n-1} &= r_n q_{n+1}.
\end{align}

Neither of the unions in (\ref{WE10WE9}) or (\ref{WE9WE10}) is
disjoint, as is evident from the proof of the lemma, where the precise value of $r_n$ did not enter in the calculation. Therefore, all $d$ and $c$ that are related by only changing the unit $r_n$ give rise to the same image and therefore to the same coset.
In appendix \ref{app:2nrep} we show that conversely, if two pairs $(c,d)$ and $(c',d')$ give rise to the same coset, then they must be related in this way. Hence, the union can be made disjoint by identifying those pairs of octavians that differ only by changing $r_n$. In the associative case this is precisely the content of Lemma~\ref{hurwlemma}.

The reason why have included $u_{r_n}$ in the definition of $w_{a,c}$, although it is an element of $W^+({\rm E}_{8})
\subset W^+({\rm E}_{9})$ and thus acts trivially on $-\delta$, comes from restricting to the associative case.
For $D_4$ and the associated Hurwitz numbers one can write $w_{a,c}$ in matrix form as
\begin{align}
w_{a,c} : X \mapsto S_{a,c}XS_{a,c}^\dagger
\end{align}
where 
$S_{a,c}$ is the ${\rm PSL}(2,\Hh)$ element
\begin{align}
S_{a,c}= T_{q_1}S_{-1}\cdots T_{q_{n+1}}S_{-1}U_{r_n}.
\end{align}
This matrix will then have the form
\begin{align}\label{ac-asterisk}
S_{a,c}=
\begin{pmatrix} 
a & b\\
c & d
\end{pmatrix}
\end{align}
where $b$ and $d$ depend on the precise choice of Euclidean decomposition.
The construction of $S_{a,c}$ shows that the left coprimality condition on $(a,c)$ is not only necessary for the existence of such a matrix in ${\rm PSL}(2,\Hh)$, but also sufficient, as we claimed in section \ref{sec:evenweyl}. The same is true for the 
left coprimality condition on $(b,d)$, and for the right coprimality conditions on $(a,b)$ and $(c,d)$.

\section{Automorphic functions on $\cH(\Div)$}
\label{sec:quatauto}

In this section, we construct Maass wave forms and Poincar\'e series
on $\cH(\Div)$ invariant under the Weyl groups studied in the preceding sections. 
We will mostly concentrate on the case $\Div=\oct$ because it is the most
interesting, and because the extension (or rather, specialization) to the other
division algebras is straightforward. The resulting expressions are simple, 
reminiscent of known expressions for ${\rm PSL}(2,\ints)$ (see for example~\cite{Goldfeld}) and reflect the arithmetic structure of the integer domains associated with the underlying Kac-Moody algebra. Our analysis can be seen as a first step towards developing a more general theory of automorphic functions and Maass wave forms of these Weyl groups.

\subsection{Maass wave forms}

The definitions of this subsection apply to all the normed division algebras (including
$\Div = \oct$). Following \cite{Iwaniec,Goldfeld} we define the scalar 
product between two functions $f,g:\,\cH(\Div)\to\cx$ by means of the invariant measure (\ref{Volume})
\begin{align}\label{L2}
(f,g) := \int_{\cal{F}} \overline{f(z)}\, g(z) \, d{\rm vol} (z).
\end{align}
A {\em Maass wave form of type $s$} with respect to the modular group $\Gamma$ 
acting on $\cH(\Div)$ is then defined to be a non-zero complex function 
$f\in{\rm L}^2({\cal{F}},\cx)$ on ${\cal{F}}\subset \cH(\Div)$ obeying
\begin{itemize}
\item $f(z) = f(\gamma(z))$ for $\gamma\in\Gamma$,
\item $  \bigtriangleup_{\text{LB}} f(z) = -s(n-s) f(z)$,
\end{itemize}
where $n={\rm dim}_\reals\,\Div$ as always, and $s\in\cx$.
If we also demand
\begin{itemize}
\item $\int_K f(u+iv) \, d^nu = 0$ for all $v>0$,
\end{itemize}
where $K$ is the `projection' of $\cal{F}$ onto $\Div$, cf. Figure~\ref{fundfig},
we would obtain so-called `cusp forms', i.e.~automorphic functions
that vanish at the cusp at infinity. However, it will be sufficient
for our purposes here to require the functions to be Fourier 
expandable around $v=\infty$.

Among the Maass wave forms one can introduce a further distinction
according to whether they are even or odd. Because $\Gamma\equiv
W_{\text{hyp}}^+$ is an index 2 subgroup in the full hyperbolic Weyl 
group $W_{\text{hyp}}$ we have two possibilities for extending the 
group action, namely 
\be
f(w(z))\, = \, \left\{ \begin{array}{c} f(z) \\[1mm]
                   \det (w) \, f(z) \end{array} \right. \;\; , 
\qquad w\in W_{\text{hyp}}\,.
\ee 
Taking the special reflection $w_\theta\in W_{\text{hyp}}$ as a reference
(see (\ref{wtheta})), we see that these definitions respectively are 
equivalent to
\be
f(-\Bz) \, = \, \left\{ \begin{array}{c} + \, f(z) \\[1mm]
                    - \, f(z) \end{array} \right.
\ee 
Maass wave forms obeying the first (second) condition are referred to 
as {\em even (odd)} Maass wave forms. Employing the generating
Weyl reflections (\ref{modularPGL}) it is straightforward to see that
odd Maass wave forms vanish on $\partial{\cal{F}}_0$, hence obey 
Dirichlet boundary conditions on ${\cal{F}}_0$ (that is, on all faces 
of the fundamental Weyl chamber), while even Maass wave forms have
vanishing normal derivatives on  $\partial{\cal{F}}_0$, hence obey
Neumann boundary conditions. If one interprets $f$ as a quantum
mechanical wave function, both possibilities are compatible with the 
rules of quantum mechanics because only scalar products of wave 
functions (for which the sign factors cancel) are observable
\cite{Kleinschmidt:2009cv,Kleinschmidt:2009hv,Forte}.

A main goal of the theory is to analyse the spectral decomposition
of the Laplace-Beltrami operator on ${\rm L}^2({\cal{F}},\cx)$
\cite{Iwaniec,Goldfeld}. As it turns out, for the discrete eigenvalues
the associated eigenfunctions must be determined numerically (see
e.g.~\cite{Then} for results on $\Div=\cx$). For odd Maass wave forms
the spectrum is purely discrete,
and one can establish the bound \cite{Kleinschmidt:2009cv,Kleinschmidt:2009hv}\footnote{In fact,
   for $SL(2,\ints)$ Maass wave forms a better bound is $E\geq 3\pi^2/2$ 
   (see \cite{Goldfeld}, p.~71), and we expect such improved lower
   bounds also to exist for the other normed division algebras $\Div$.}
\be
- \bigtriangleup_{\text{LB}} \geq \frac{n^2}4
\ee
generalizing an argument of \cite{Iwaniec}. For even $f(z)$ 
one has in addition to the discrete spectrum a continuous spectrum that
can be constructed by means of Eisenstein series and extends from
$n^2/4$ over the positive real axis (see below).

\subsection{Poincar\'e series}
\label{sec:eisen}

We recall that the Laplacian on hyperbolic space $\cH(\Div)$ of dimension 
$n+1$ is given by (\ref{LB}), again with $n= {\rm dim}_\reals\, \Div$.
For $I_s(z)= v^s$ or $I_{n-s}(z) = v^{n-s}$ with $s\in\cx$ we compute
\begin{align}
- \bigtriangleup_{\text{LB}} \, I_s(z) = E_s\, I_s(z) 
     = s(n-s) \, I_s(z).
\end{align}
Note that for real values of $s$  we have $E_s>0$ if and only if $0<s<n$.

We are here interested in functions that are both eigenfunctions of
the Laplacian $\bigtriangleup_{\text{LB}}$ {\em and} automorphic
(but not necessarily square integrable). To turn the 
function $I_s$ into an automorphic form, we have to make it symmetric 
under modular transformation by `averaging' over the modular group
$\Gamma\equiv W_{\text{hyp}}^+$. Abstractly, this is achieved by 
defining the {\em Poincar\'e series} (or restricted Eisenstein series) 
\begin{align}\label{poin}
\cP_s(z) : = \sum_{\gamma\,\in \, \Gamma_\infty\backslash \Gamma} 
I_s(\gamma(z)).
\end{align}
Whenever this sum converges it is an eigenfunction of the 
Laplace-Beltrami  operator by virtue of the invariance property 
(\ref{LB1}) which implies that every term in the sum is separately 
an eigenfunction with the same eigenvalue. The restriction to summing 
over cosets of the stabilizer of $I_s(z)$ ensures that the sum is 
well-defined. Furthermore $I_s(w_\theta(z)) 
= I_s( -\Bz) = + I_s(z)$, so we conclude that $\cP_s(z)$ 
satisfies Neumann boundary conditions on ${\cal{F}}_0$.

We emphasize again that the above definition, 
though standard \cite{Iwaniec,Goldfeld}, is special inasmuch
all our  `modular groups' $\Gamma$ are hyperbolic Weyl groups,
with the corresponding affine Weyl subgroups $\Gamma_\infty$ as the 
stabilizer groups of the cusp at infinity. The relation of the hyperbolic Weyl group to integer domains in $\Div$ will bring the arithmetic structure to the fore in the final expressions for $\cP_s$. We will perform the analysis for $\Gamma=W^+({\rm E}_{10})$ since for example the quaternionic case $W^+({\rm D}_4^{++})$ is contained in it by specialization.

In order to evaluate the sum in (\ref{poin}) we employ 
the description (\ref{WE9WE10})
of the
right cosets $\Gamma_\infty\backslash \Gamma$. 
Consider therefore two left coprime octavians $c$ and $d$ with corresponding left Euclidean decomposition (\ref{leuclalg}) and associated left coset representative 
\begin{align}
\tilde{w}_{c,d} = u_{\bar{r}_n} \circ s_{-1}\circ t_{q_{n+1}}  \circ\dots\circ s_{-1}\circ t_{q_1}.
\end{align}
as in (\ref{wcd}).

\begin{lemma}
Let $(z_i)$ be a sequence of elements in $\cH(\Div)$ defined recursively by $z_0=u_0+i v_0:= z$ and $z_{k+1} =u_{k+1} +i v_{k+1} := (s_{-1}\circ t_{q_{k+1}})(z_k) = -(z_k+q_{k+1})^{-1}$. 
Then
\be
v_{n+1} = \frac{v_i}{|r_i z_i+r_{i-1}|^2}
\ee
for all $i=n+1$ down to $i=0$. (Here, we let $r_{n+1}=0$, $r_0=c$ and $r_{-1}=d$ as before.)
\end{lemma}

\Pf For $i=n+1$ the claim is true since $r_{n+1}=0$ and $|r_n|=1$. The induction step consists of
\begin{align}
v_{n+1} &= \frac{v_i}{|r_i z_i +r_{i-1}|^2} =  \frac{v_{i-1}}{ |r_i z_i +r_{i-1}|^2 |z_{i-1}+q_i|^2}\nn\\
&= \frac{v_i}{ |-r_i(z_{i-1}+q_i)^{-1}+r_{i-1}|^2 |z_{i-1}+q_i|^2}\nn\\
&= \frac{v_{i-1}}{|r_{i-1} z_{i-1}+r_{i-1}q_i-r_i|^2}\\
&= \frac{v_{i-1}}{|r_{i-1}z_{i-1}+r_{i-2}|^2}.\nn
\end{align}
In one step one requires that $r_{i-1}=[r_{i-1}(z_{i-1}+q_i)](z_{i-1}+q_i)^{-1}$. 
Remarkably, this step remains valid even for $\cH(\oct)$: though such a relation is 
not true generally over the sedenions due to lack of alternativity, it holds in the present 
case since the imaginary part of $z_{i-1}+q_i$ is $i$ multiplied with
a real number only. 
For the same reason one can also use multiplicativity of the norm.\qed

\noindent
Since the imaginary part of $z_{n+1}$ is not changed by the final $u_{r_n}$ we obtain
\begin{align}
I_s\left(\tilde{w}_{c,d}(z)\right)= v_{n+1}^s  = \frac{v_0^s}{|r_0 z_0 + r_{-1}|^{2s}} = \frac{v^s}{|cz+d|^{2s}}\,.
\end{align}
The fact that $r_n$ does not change the result of this computation 
indicates that one can choose the unit $r_n$ freely. 
As we stated in section \ref{we9we10-section}, if and only if $c'$ and $d'$ are defined by the same $q_i$ as $c$ and $d$ but with $r_n'$ instead of $r_n$ the corresponding $\tilde{w}_{c',d'}$ will be in the same coset as $\tilde{w}_{c,d}$ and hence $I_s(\tilde{w}_{c',d'}(z))=I_s(\tilde{w}_{c,d}(z))$. Therefore we can undo the overcounting associated with the non-disjointness of (\ref{WE9WE10}) by dividing by the number of choices for $r_n$, that is, the number of units of the integer domain in $\Div$.
Finally, we obtain the following explicit expression for the Poincar\'e series
\begin{align}\label{eisen}
\cP_s(z) = \frac{1}N
\sum_{c,d\in\cO \,\,\text{left coprime}} \frac{v^s}{|cz+d|^{2s}}\,.
\end{align}
where $N$ is the number of units, which is $240$ for the octavians and $24$ for the Hurwitz numbers.
This expression crucially uses the number theoretic properties of the integer domain $\cO\subset\Div$ underlying the hyperbolic Weyl group.

One can also define the {\em unrestricted} Eisenstein series
\begin{align}\label{QEisen}
\cE_s(z)
:= {\sum_{(c,d)\in\cO^2\backslash \{(0,0)\}}} \frac{v^s}{|cz + d|^{2s}} \,.
\end{align}
This  Eisenstein series is related to (\ref{eisen}) by
\begin{align}
\cE_s(z) = \zeta_{\cO}(s) \cP_s(z)\,,
\end{align}
where 
\begin{align}
\zeta_{\cO}(s) = \sum_{0\neq a\in\cO} |a|^{-s}
                     = \sum_{n\in\mathbb{N}} \frac{\sigma_\cO(n)}{n^s}
\end{align}
is the (Dedekind) zeta function associated with the appropriate integers
and $\sigma_\cO(n)$ counts the number of roots of (squared) length 
$2n$ on the root lattice (which is always a multiple of the number of units). This zeta function is also the correct factor of proportionality in the non-associative case as can be seen by considering any pair $(c,d)$. The last vanishing remainder $r_n$ is not necessarily a unit but has well-defined norm $|r_n|^2$. Then one can again construct all pairs $(c',d')$ which yield the same summand in (\ref{QEisen}) by letting $r_n$ range over all octavians of given norm $|r_n|^2$. This produces exactly the zeta function (encoding the same information as the ${\rm E}_8$ lattice theta function).

We expect that with an appropriately defined completed zeta
function ${\xi}_{\cO}(s)$ the following functional relation 
holds\footnote{For example, the completed Riemann zeta function is
\begin{align}
\xi(s) = \pi^{-s/2} \Gamma(s/2) \zeta(s)\nn
\end{align}
and satisfies the functional relation $\xi(s)=\xi(1-s)$ 
which can be used to define the Riemann zeta function by analytic 
continuation outside its domain of convergence $\text{Re}(s)>1$. 
There exist appropriate generalizations for other algebraic integers.}
\begin{align}\label{funceq}
\xi_{\cO}(s) \cP_s(z) = \xi_{\cO}(n-s) \cP_{n-s}(z)\,,
\end{align}
where $n$ is the dimension of the division algebra and this functional relation should be related to the ones studied in the real rank one case in~\cite{Langlands}.

We now turn to the Fourier expansion of the Poincar\'e series. Because the sum (\ref{poin}) is periodic under integer shifts $u\to u+o$ for $o\in\cO$, it can be expanded into a Fourier series. Since it is an eigenfunction of the Laplacian, the Fourier coefficients have to also satisfy differential equations and we obtain
\begin{align}\label{fourexp}
\cP_s(z) = v^s + a(s) v^{n-s}+   
v^{n/2} \sum_{\mu\in \cO^*\backslash\{0\}} a_\mu K_{s-n/2} 
(2\pi|\mu|v) e^{2\pi i \mu(u)}\,,
\end{align}
where $\cO^*\subset \Div$ is the lattice dual to the lattice of 
integers $\cO\subset\Div$ relevant to the hyperbolic Weyl group 
and $K_\nu(v)$ is the solution to the Bessel equation 
\begin{align}
v^2 \,\frac{d^2K_\nu}{dv^2} + v \, \frac{dK_\nu}{dv} + (v^2-\nu^2)\, K_\nu = 0,
\end{align}
which vanishes exponentially for large $v$. The coefficients 
$a_\mu$ in (\ref{fourexp}) are further constrained by the even part of the finite Weyl 
group which acts as a set of (generalized) rotations on $\cO$, so 
that $a_\mu= a_{s(\mu)}$ for all $s\in W^+_{\text{fin}}$. Further 
constraints come from a Hecke algebra (if it can be suitably defined) 
and will render the coefficients $a_\mu$ multiplicative over 
$\Div=\reals,\cx$. The first two terms in (\ref{fourexp}) correspond 
to the so-called constant terms. For equation (\ref{funceq}) to hold one needs $a(s)=\xi_{\cO}(n-s)/\xi_{\cO}(s)$. The remaining terms fall of
exponentially as $v\rightarrow\infty$. In string theory applications, $v$ is associated with the string dilaton and hence  the first two terms correspond to perturbative effects whereas the exponentially suppressed terms are non-perturbative in the string coupling.

By standard arguments one can show that (\ref{eisen}) and (\ref{QEisen}) are convergent for $\text{Re}(s)> n/2$. However, these functions are never square integrable with respect to (\ref{L2}). 
For both ${\rm Re}\, s >n/2$ and ${\rm Re}\, s< n/2$, 
the functions $\cE_s(z)$  are not  normalizable with respect to the 
invariant measure (\ref{Volume}) when integrated over the fundamental 
domain as in (\ref{L2}). This is due to the divergence of the $v$-integral 
for $v\rightarrow\infty$ because of the `constant terms' $v^s$ 
and $v^{n-s}$, both of which appear with non-vanishing coefficients 
in (\ref{fourexp}), whereas the `tail'  involving Bessel functions decays 
exponentially for large $v$. However, for the special values on the `critical line'
($r\in\reals$)
\begin{align}
s= \frac{n}2 + ir
\end{align}
we obtain `almost' normalizable states with eigenvalues
\begin{align}
E= \frac{n^2}4 + r^2 \geq \frac{n^2}4.
\end{align}
Notwithstanding a rigorous proof, this can be
roughly seen as follows. Exploiting the expansion (\ref{fourexp})
the leading terms behave like $v^s$ or $v^{n-s}$, with the subleading 
terms being integrable in $v$. Substituting these `dangerous' 
terms into (\ref{L2}) and neglecting finite contributions, 
we get (with ${\rm Re}\, s= n/2$)
\begin{align}
\int_{\cal{F}}  d{\rm vol}(z) \, \overline{\cE_s (z)} \cE_{s'}(z) 
\sim \int_K d^n u \int^\infty_{\sqrt{1-|u|^2}}  \frac{dv}{v} \exp\big[i(\pm r' \pm r) \ln v\big].
\end{align}
The integral over $u$ is finite because $K$ is compact. Changing 
variables to $\xi = \ln v$ and extending the range of integration
to $v=0$, we get
\begin{align}
\int^\infty_{-\infty} d\xi\, 
\exp\big[i(\pm r' \pm r)\xi\big] \; = \; 2\pi \delta(r' \pm r).
\end{align}
Hence, for these values of $s$, the Eisenstein functions are
$\delta$-function normalizable up to finite corrections.
As a result, there is a continuous part of the spectrum
\begin{align}
{\rm spec}_c ( - \bigtriangleup_{\text{LB}}) =
  \left[ \frac{n^2}4 \, ,\, \infty \right)
\end{align}
for {\em even} Maass wave forms. In addition there are discrete
eigenvalues for even Maass wave forms, which are embedded
in the continuum (but which must be determined numerically).
The spectrum of odd Maass wave forms is purely discrete.

\subsubsection*{Acknowledgments}
We would like to thank Jonas T.~Hartwig and Michael Koehn for discussions.
AK is a Research Associate of the Fonds de la Recherche Scientifique--FNRS, Belgium.
The work by AK and JP is in part supported by IISN - Belgium (conventions 4.4511.06 and 4.4514.08), by the Belgian Federal Science Policy Office 
through the Interuniversity Attraction Pole P6/11. AK and JP would also like to thank the Albert Einstein Institute for hospitality.

\appendix

\section{Automorphisms of octonions} \label{autoctapp}
In this appendix, we will prove a theorem about the automorphisms of the octonions that we use in 
section~\ref{evene7weylgroup}.
First we list the Moufang identities that we also use
\begin{align}
(ax)(ya)&=a(xy)a,&
((xa)y)a)&=x(aya),&
a(x(ay))&=(axa)y. \label{moufang}
\end{align}
These identities follow from the alternative laws (\ref{altlaws}). (For a proof, see \cite{schafer}.)
Another useful identity that follows from the alternative laws and the Moufang identities is
\begin{align}
(a^2 x a)(a^{-1} ya)=a^2(xy)a.
\end{align}
Indeed, we have 
\begin{align} \label{alt-identity}
(a^2 x a)(a^{-1} ya)&=(a (ax) a)(a^{-1} ya)\nn\\
&=a ((ax a)(a^{-1} y))a\nn\\
&=a (a(xy))a=a^2(xy)a.
\end{align}
Now we can prove the following lemma.
\begin{lemma}\label{lemma1} For any integer $n\leq 1$ and any $z,\,a_1,\,a_2,\,\ldots,\,a_n \in \oct$ we have
\begin{align}
&a_1 (a_2 ( \cdots (a_n z a_n^{-1}) \cdots )a_2^{-1})a_1^{-1}\nn\\
&\quad=
a_1^{-2}(a_2^{-2}\cdots(a_n^{-2}(b_nz)a_n^{-1})\cdots a_2^{-1})a_1^{-1} \label{ablemma}
\end{align} 
where
\begin{align}
b_n&=a_n^2(a_{n-1}^2\cdots(a_2^2 a_1^3 a_2)\cdots a_{n-1})a_n.
\end{align}
\end{lemma}
\Pf We prove this by induction over $n$. The case $n=1$ follows easily by alternativity. Suppose now that 
the identity (\ref{ablemma}) holds for some integer $n\geq 1$. Then we have
\begin{align}
&a_1 (a_2 \cdots (a_n (a_{n+1} z a_{n+1}^{-1})a_n^{-1}) \cdots a_2^{-1})a_1^{-1}\nn\\
&\quad=
a_1^{-2}(a_2^{-2}\cdots(a_n^{-2}(b_n(a_{n+1} z a_{n+1}^{-1}))a_n^{-1})\cdots a_2^{-1})a_1^{-1} 
\end{align} 
so we only have to show that
\begin{align} \label{ablemmaproof}
b_n(a_{n+1}za_{n+1}^{-1})=a_{n+1}^{-2}(b_{n+1}z)a_{n+1}^{-1}.
\end{align}
To do this we use that
\begin{align}
b_n = a_{n+1}^{-2}b_{n+1}a_{n+1}^{-1}
\end{align}
so that the left hand side of (\ref{ablemmaproof}) becomes
\begin{align}
b_n(a_{n+1}za_{n+1}^{-1})=
(a_{n+1}^{-2}b_{n+1}a_{n+1}^{-1})(a_{n+1}za_{n+1}^{-1}).
\end{align}
and then (\ref{ablemmaproof}) follows from (\ref{alt-identity}).
\qed
\begin{lemma}\label{lemma2} For any integer $n\leq 1$ and any $x,\,y,\,a_1,\,a_2,\,\ldots,\,a_n \in \oct$ we have
\begin{align} 
&(a_1 ( \cdots (a_n x a_n^{-1}) \cdots )a_1^{-1})
(a_1 ( \cdots (a_n y a_n^{-1}) \cdots )a_1^{-1})\nn\\
&\quad=
a_1^{-2}(a_2^{-2}\cdots(a_n^{-2}((b_nx)y)a_n^{-1})\cdots a_2^{-1})a_1^{-1} \label{lemma2-eq}
\end{align} 
where
\begin{align}
b_n&=a_n^2(a_{n-1}^2\cdots(a_2^2a_1^3a_2)\cdots a_{n-1})a_n.
\end{align}
\end{lemma}
\Pf By Lemma \ref{lemma1} we have
\begin{align}
&(a_1 ( \cdots (a_n x a_n^{-1}) \cdots )a_1^{-1})
(a_1 ( \cdots (a_n y a_n^{-1}) \cdots )a_1^{-1})\nn\\
&\quad= (a_1^{-2}\cdots(a_n^{-2}(b_nx)a_n^{-1})\cdots a_1^{-1})(a_1 ( \cdots (a_n y a_n^{-1}) \cdots )a_1^{-1})
\end{align}
and then (\ref{lemma2-eq}) follows by using (\ref{alt-identity}) successively.
\qed\\
\noindent
No we can prove the main result of this appendix.
\begin{thm} \label{thetheorem}
A transformation
\begin{align}
\varphi: x \mapsto a_1 (a_2 ( \cdots (a_n x a_n^{-1}) \cdots )a_2^{-1})a_1^{-1}
\end{align}
of the octonions, where $a_1,\,\ldots,\,a_n\neq0$, is an automorphism if and only if
\begin{align}
b_n = a_n^{2}(a_{n-1}^{2}(\cdots (a_{2}^2 a_1^3 a_{2}) \cdots a_{n-1})a_n \in \reals.
\end{align}
\end{thm}
\Pf Setting $xy=z$, we have in Lemma \ref{lemma1} and Lemma \ref{lemma2} obtained expressions for $\varphi(xy)$ and $\varphi(x)\varphi(y)$, respectively. The difference between these two expressions is
\begin{align}
&\varphi(xy)-\varphi(x)\varphi(y)
= a_1^{-2}(a_2^{-2}\cdots(a_n^{-2}\{b_n,\,x,\,y\}a_n^{-1})\cdots a_2^{-1})a_1^{-1}.
\end{align}
Since there are no zero divisors, this difference is zero if and only if the associator $\{b_n,\,x,\,y\}=b_n(xy)-(b_n x)y$ vanishes. Since this must happen for all $x$ and $y$, the necessary and sufficient condition is $b_n \in \reals$.
\qed\\
\noindent
The theorem simplifies considerably if we take $a_1,\,\ldots,\,a_n$ to be unit octonions, and furthermore imaginary.
\begin{cor}\label{thecorollary}
A transformation
\begin{align}
\varphi: x \mapsto a_1 (a_2 ( \cdots (a_n x a_n) \cdots )a_2)a_1
\end{align}
of the octonions, where $a_1,\,\ldots,\,a_n$ are imaginary unit octonions or real
unit  octonions, is an automorphism if and only if
\begin{align}
(((a_1 a_2)a_3)\cdots a_n) = \pm 1.
\end{align}
\end{cor}
\Pf This follows directly from Theorem \ref{thetheorem}
if we use that
$a^{-1}=\pm a$
and $a^2=a^{-2}=\pm 1$
for any imaginary unit octonion or real unit octonion $a$.
\qed\\
\noindent
We use this corollary in section \ref{evene7weylgroup} when we identify the imaginary unit octavians with the simple roots of ${\rm E}_7$, and the corresponding bimultiplications with generators of $W^+({\rm E}_7)$.

\section{$2n$-dimensional representation of $W^+_{\text{hyp}}$} \label{app:2nrep}

The representation (\ref{evenweylgen}) of the even hyperbolic Weyl group as conjugation on Hermitian $(2\times 2)$-matrices over $\Div$ is an $(n+2)$-dimensional representation, corresponding to the action on the Cartan subalgebra. Alternatively, it is the action in the vector (fundamental) representation of ${\rm SO}(1,n+1)$ in which $W^+_{\text{hyp}}$ is embedded as a discrete subgroup. The isomorphism $W^+_{\text{hyp}}\cong {\rm PSL}(2,\cO)$ for integers $\cO\subset\Div$ suggests also a natural action on $2$-component vectors of elements in $\Div$. Such a representation corresponds to the spinor representation of ${\rm SO}(1,n+1)$.

Concretely, we consider the {\em right} action of ${\rm PSL}(2,\cO)$ on the $2n$-dimensional space of two components {\em rows}
\be
\left\{ (a_1, a_2)\,:\, a_1,a_2 \in \Div \right\}/\sim\,,
\ee
where $\sim$ denotes the equivalence relation associated with the `${\rm P}$' in ${\rm PSL}(2,\cO)$. Except for the complex case it corresponds to $(a_1,a_2)\sim (-a_1,-a_2)$. The action of ${\rm PSL}(2,\cO)$ is by right multiplication of the generators (\ref{evenweylgen}).

In order to check that this defines a representation of $W^+_{\text{hyp}}$ one has to verify the defining relations. We do this in a non-associative example for the product $S_i S_j$ of two generators such that $(\ve_i\ve_j)^2=\pm 1$. From the definition one gets
\be
(a_1,a_2)\cdot (S_iS_j)^2= \big( (((a_1\ve_i)\ve_j)\ve_i)\ve_j\,,(((a_2\bar{\ve}_i)\bar{\ve}_j)\bar{\ve}_i)\bar{\ve}_j\big)\,.
\ee
Now, using the Moufang identities one calculates 
\be
(((a_1\ve_i)\ve_j)\ve_i)\ve_j\ &=& (a_1\ve_i)(\ve_j\ve_i\ve_j) = (a_1(\bar{\ve}_j\bar{\ve}_i\bar{\ve}_j))(\ve_j\ve_i\ve_j\ve_i\ve_j\ve_i\ve_j)\nn\\
&=&\pm (a_1(\bar{\ve}_j\bar{\ve}_i\bar{\ve}_j) (\ve_j\ve_i\ve_j) = \pm a_1\,,
\ee
where alternativity has been used in a number of places. Hence, we arrive at $(S_iS_j)^2={\text{id}}\in {\rm PSL}(2,\cO)$ as required. The calculations in the other cases are similar.

One advantage of this representation is that one can easily exhibit the invariant content of elements like $\tilde{w}_{c,d}$ in (\ref{wcd}). A calculation shows that
\be
(0,1)\cdot \tilde{w}_{c,d} = (c,d)\,.
\ee
In the associative case this is clear since $\tilde{w}_{c,d}$ is then given by a matrix with 
bottom row $(c,d)$. In the non-associative case one cannot write a single matrix for $\tilde{w}_{c,d}$ but the above computation shows that it still behaves as a matrix with bottom row $(c,d)$ when it acts on $(0,1)$.\footnote{In this representation, the stabilizer of $(0,1)$ is only the translation group $\mathcal{T}$ and not the full (even) affine Weyl group $\Gamma_\infty\equiv W^+_{\text{aff}}=W_{\text{fin}}^+\ltimes \mathcal{T}$. Therefore we can `resolve' the presence of $u_{r_n}$ in (\ref{wcd}).}
This insight allows us to derive the overcounting factor $N$ appearing in (\ref{eisen}).
Suppose that $(c',d')$ is a pair of left coprime octavians that gives rise to the same coset as $(c,d)$, that is,
\begin{align}
W^+({\rm E}_9)\tilde{w}_{c,d} = W^+({\rm E}_9)\tilde{w}_{c',d'}
\end{align}
Since
$\Gamma_\infty=W^+_{\text{aff}}$ is generated (non-minimally) by the translations $t_y$ of (\ref{transdef}) and the rotations $u_\ve$ of (\ref{rotgen}) this means that $\tilde{w}_{c',d'}=w\tilde{w}_{c,d}$, where $w$ is a product of such translations  
and rotations. Acting on $(0,1)$ we get
\begin{align}
(c',d')=(0,1)\cdot (w\tilde{w}_{c,d}).
\end{align}
The rotations in $w$ will change $(0,1)$ to $(0,\varepsilon)$ where $\varepsilon$ is a unit, and the translations act trivially on such rows.

We conclude that
$c'$ and $d'$ are determined by the left Euclidean algorithm with the same $q_i$ as for $(c,d)$ but with $r_n'=\ve r_n$ instead of $r_n$. Therefore all pairs that represent the same coset must be related in such a way. One can show that for all choices of $\ve$ one obtains distinct pairs, and therefore the overcounting in (\ref{eisen}) is equal to the number of units in $\cO$.

\section{Green functions on $\cH(\Div)$}
\label{app:Green}

For $z_1, z_2 \in \cH(\Div)$ define (cf. (\ref{Length}))
\be
\la(z_1,z_2) := \frac{|u_1 - u_2|^2 + (v_1 -v_2)^2}{4v_1 v_2}
        \equiv \frac12 \sinh^2 \frac{d(z_1,z_2)}{2}
\ee
Acting with the Laplace-Beltrami operator (\ref{LB}) on the first 
argument of $\la(z_1,z_2)$  we get
\be
\LB \, \la(z_1, z_2) = (n+1) \left( \frac12 + \la(z_1,z_2) \right)
\ee
and thus for $\xi\in\reals$ and $\la\equiv\la(z_1,z_2)$
\be
\LB \, \big[\xi + \la\big]^{-s} &=& s\big[\xi + \la\big]^{-s-2} \times \nn\\
&&\!\!\!\! \!\!\!\! \!\!\!\! \!\!\!\! \!\!\!\! \!\!\!\! \!\!\!\! \!\!\!\! \!\!\!\! \!\!\!\!  \times \;
\left[ (s+1)(\la + \la^2)  
- (n+1)\left( \frac12 + \la\right) (\xi + \la)\right]
\ee
Following \cite{Iwaniec}, we define the Green 
function (or the `propagator') as
\be
G_s(\la(z,w)) &=& \int_0^1 \left[ \xi(1-\xi)\right]^{s- \frac{n+1}2}
                  \big[\xi + \la(z,w)\big]^{-s} d\xi 
\ee
for $z\,,\,w\in\cH(\Div)$. This integral converges for $s > \frac{n-1}2$
and $\la(z,w) >0$.  We will assume $n>1$ from now on, as the case
$n=1$ is treated in much detail in \cite{Iwaniec}.

For non-coincident points $z\neq w$, application
of the operator $\LB + s(n-s)$ gives an integrand proportional to
\be
s\frac{d}{d\xi}\left[\Big(\xi(1-\xi)\Big)^{s+\frac{1-n}2}
           \big(\xi + \la(z,w)\big)^{-s-1}\right]
\ee 
Therefore the boundary terms, and hence the integral vanish for 
$s > \frac{n-1}2$. This shows that, for $z\neq w$, we have
\be\label{LBonG}
\Big[ \LB + s(n-s) \Big] G_s\big(\la(z,w)\big) = 0
\ee
For coincident arguments the integrand is singular, and we must
reason more carefully. To determine the behavior for small $\la$
we split the integral as
\be
\int_0^1 = \int_0^A + \int_A^B + \int_B^1
\ee
with suitable $0 < A < B < 1$. Since for $\la\geq 0$ and $0<\xi <1$
\be
\left| \frac{\xi(1-\xi)}{\xi + \la}\right| < 1
\ee
the first integral is bounded above by
\be
\int_0^A \frac{d\xi}{(\xi + \la)^{\frac{n+1}2}}
  = A \la^{- \frac{n+1}2} + \cO \big( A^2\la^{-\frac{n+3}2}\big)
\ee
Choosing $A=\la^{\frac{n+1}2}$, we see that this part of the integral
is $\cO(1)$ for small $\la$. For the middle integral we have 
\be
\int_A^B   \frac{d\xi}{(\xi + \la)^{\frac{n+1}2}}
= \frac2{n-1} \left[ \frac1{(\la + A)^{\frac{n-1}2}} 
            - \frac1{(\la + B)^{\frac{n-1}2}}\right] \,.
\ee
so that with (say) $B=\frac12$ we get
\be\label{singular}
\int_{\la^{(n+1)/2}}^{1/2}    \frac{d\xi}{(\xi + \la)^{\frac{n+1}2}} =
\frac2{n-1} \la^{-\frac{n-1}2} + \cO(1) \quad
\mbox{for $\la\rightarrow 0$.}
\ee
The remaining integral is again bounded as
\be
\int_{1/2}^1 (1-\xi)^{s- \frac{n+1}2} d\xi = \cO(1)
\ee
for $s> \frac{n-1}2$. From (\ref{singular}) we see that the
Green function behaves as $\sim \la^{-(n-1)/2}$ for small
$\la$, and therefore eq.~(\ref{LBonG}) must be amended to
\be
\Big[ \LB + s(n-s) \Big] G_s\big(\la(z,w)\big) = 
\frac{4\pi^{\frac{n+1}2}}{(n-1)\Gamma\left(\frac{n+1}2\right)}  \, \delta^{(n+1)}(z,w)
\ee
where we have used the standard formula for the volume of 
the $(n+1)$-dimensional unit sphere.
The automorphic Green function is then defined as
\be
\mathbb{G}_s(z/w) := \sum_{\gamma\in\Gamma} G_s(z, \gamma(w))
\ee
in analogy with the case $n=1$ \cite{Iwaniec}.

\section{Geodesics in $\cH(\quat)$ and periodic orbits}
\label{app:geo}

Geodesics in $\cH(\quat)$ are given by half-circles and straight
lines. For the usual complex upper half plane it
is known (see e.g. \cite{Bogo}) that we can associate to each
hyperbolic element\footnote{Recall that a 
  hyperbolic element  $M\in {\rm SL}(2,\reals)$ has reciprocal real 
  eigenvalues or, equivalently, satisfies $|\rm{Tr}\,M|>2 $.} 
$M\in {\rm PSL}(2,\reals)$ a geodesic that is mapped onto
itself by the action of $M$ (the same is true for all matrices
conjugate to $M$). 
Let us consider the
`imaginary axis' in $\cH(\quat)$, parametrized by $z(t) = it$ with
$0<t<\infty$. For this geodesic, the hyperbolic motion leaving it
invariant is obviously given by $\gamma_t = {\rm diag}\, (t^{1/2} ,
t^{-1/2})$, whereby $i$ is mapped to $\gamma_t(i)= \gamma_t (z(1))=
it$. It is straightforward to see that all geodesics in $\cH(\quat)$ 
are ${\rm PSL}(2,\quat)$ images of one another. First,
straight line geodesics centered at $u\neq0$ are trivially 
obtained by acting on $z(t) = it$ with the shift matrix
\be\label{Su}
T_u = \begin{pmatrix}
1 & u\\
0 & 1
\end{pmatrix}\,.
\ee
For arbitrary $u_1 \neq u_2\in\quat$ the circular geodesic with endpoints 
$u_1$ and $u_2$ can (for example) be obtained by acting on 
$z(t)$ with the matrix
\be\label{Su1u2}
C_{u_1,u_2} = \frac1{\sqrt{|u_1-u_2|}}\begin{pmatrix}
u_2 & u_1\\
1 & 1
\end{pmatrix}\,.
\ee
Use of (\ref{z'}) results in the explicit parametrization
of the geodesic half-circle 
\begin{align}\label{geo}
z'(t) = \frac{u_1 +  u_2 \, t^2\, + \, it\, |u_1- u_2|}{1 + t^2}
\qquad (0<t<\infty) \; .
\end{align}
The hyperbolic motions leaving invariant this geodesic half-circle 
are now simply given by (for all $t$)
\be\label{Mt}
M_t = S \, \begin{pmatrix}
         t^{1/2} & 0\\
          0 & t^{-1/2}
\end{pmatrix}  \, S^{-1}\,,
\ee
where $S$ is either $T_u$ or $C_{u_1,u_2}$.
The formula expressing $M_t$ explicitly in terms of $t, u_1 , u_2$
follows directly from (\ref{Su1u2}) and (\ref{Sinverse}) but is not 
very illuminating. However, it implies the inequality
\be\label{ReMt}
{\rm Re}\,\big( {\rm Tr}\, M_t \big) = t^{1/2} + t^{-1/2} \geq 2\,.
\ee
Note that the real part of the trace is cyclic 
even over quaternionic matrices and hence furnishes an invariant for 
matrices up to conjugation. In other words, we may adopt (\ref{ReMt})
as the quaternionic generalization of the usual condition of 
hyperbolicity for an ${\rm SL}(2,\reals)$ matrix. 

Consider now the case when $M\equiv M_{t_0}$ happens to be integral
for some value $t_0 > 1$, i.e., $M_{t_0}\in {\rm PSL}^{(0)}(2,\Hh)$. Then the associated geodesic 
gives rise to a {\em periodic orbit in the fundamental domain $\cal{F}$} 
of the modular group ${\rm PSL}^{(0)}(2,\Hh)$. This is most easily seen by following 
the geodesic circle successively through the images of the fundamental 
domain $\cal{F}$ until we reach the image  of $\cal{F}$ under the special 
transformation $M$, where we connect up again to the original geodesic 
curve intersecting  $\cal{F}$. (The integrality condition 
$M\in {\rm PSL}^{(0)} (2,\,\Hh)$ may be viewed as the analog of the 
rationality condition for periodic orbits on tori.) The length $\ell_p$ 
of this periodic orbit is easily calculated by mapping it back
to the imaginary axis, whence 
\be
\ell_p = \int_1^{t_0} \frac{dv}{v} = \log t_0
\ee
by the invariance of the geodesic length under the modular group.
With (\ref{ReMt}) we recover the quaternionic analog of
the well known formula for ${\rm PSL}(2,\ints)$ \cite{Bogo}
\be
2 \cosh \frac{\ell_{p}}{2} = \big|\, \text{Re} \,(\text{Tr}\, M)\, \big|\,.
\ee
The limiting case $\ell_p\to\infty$ corresponds to the infinite geodesic 
along the imaginary axis. Consequently, the number of periodic 
orbits in the fundamental domain  $\cal{F}$ increases exponentially 
with the geodesic length of the orbit, as it does for  ${\rm PSL}(2,\ints)$ 
(see e.g. \cite{Bogo}).

\begingroup\raggedright

\endgroup

\end{document}

%% file: uhp5.pstex_t
\begin{picture}(0,0)%
\includegraphics{uhp5.pstex}%
\end{picture}%
\setlength{\unitlength}{3947sp}%
\begingroup\makeatletter\ifx\SetFigFont\undefined%
\gdef\SetFigFont#1#2#3#4#5{%
  \reset@font\fontsize{#1}{#2pt}%
  \fontfamily{#3}\fontseries{#4}\fontshape{#5}%
  \selectfont}%
\fi\endgroup%
\begin{picture}(6187,4861)(2979,-4325)
\put(5851,389){\makebox(0,0)[lb]{\smash{{\SetFigFont{12}{14.4}{\rmdefault}{\mddefault}{\updefault}{\color[rgb]{0,0,0}$v$}%
}}}}
\put(9151,-4036){\makebox(0,0)[lb]{\smash{{\SetFigFont{12}{14.4}{\rmdefault}{\mddefault}{\updefault}{\color[rgb]{0,0,0}$u\in\mathbb{A}$}%
}}}}
\put(5851,-4261){\makebox(0,0)[lb]{\smash{{\SetFigFont{12}{14.4}{\rmdefault}{\mddefault}{\updefault}{\color[rgb]{0,0,0}$K$}%
}}}}
\end{picture}%